\documentclass[opre,nonblindrev]{informs3} 

\DoubleSpacedXI 

\usepackage{endnotes}
\let\footnote=\endnote

\newcommand{\va}{\textcolor{red}}       

\usepackage{natbib}
 \bibpunct[, ]{(}{)}{,}{a}{}{,}%
 \def\bibfont{\small}%
 \def\bibsep{\smallskipamount}%
 \def\bibhang{24pt}%
\TheoremsNumberedThrough     
\ECRepeatTheorems

\EquationsNumberedThrough    

\begin{document}

\TITLE{Online Resource Allocation with Customer Choice}
\ARTICLEAUTHORS{%
\AUTHOR{Guillermo Gallego, Anran Li, Van-Anh Truong, Xinshang Wang} \AFF{Department
of Industrial Engineering and Operations Research, Columbia
University, New York, NY, USA, \EMAIL{gmg2@columbia.edu,al2942@columbia.edu, vatruong@ieor.columbia.edu,xw2230@columbia.edu}
\URL{}}}

\ABSTRACT{We introduce a general model of resource allocation with customer choice.  In this model, there are multiple resources that are available over a finite horizon. 
The resources are non-replenishable and perishable.    Each unit of a resource can be instantly made into one of several products.  There are multiple customer types arriving randomly over time.  An assortment of products must be offered to each arriving customer, depending on the type of the customer, the time of arrival, and the remaining inventory.  From this assortment, the customer selects a product according to a general choice model.  The selection generates a product-dependent and customer-type-dependent reward.  The objective of the system is to maximize the total expected reward earned over the horizon.  

The above problem has a number of applications, including personalized assortment optimization, revenue management of parallel flights, and web- and mobile-based appointment scheduling.  We derive online algorithms that are asymptotically optimal and achieve the best constant relative performance guarantees for this class of problems.}
\maketitle

\section{Introduction}

In this paper we introduce a general model of resource allocation with customer choice.  In this model, there is a finite, continuous-time horizon.  Over the horizon, there is a known set of resources, each finite in quantity and perishing at a known date.  Each unit of a resource can be used to instantly make one of several products.   Over the horizon, customers of various types arrive according to non-homogenous Poisson processes.  The type of a customer is observable by the system.  Upon a customer's arrival, according to the time of arrival, the inventory of available resources, and the type of the customer, the system chooses an assortment of products to display.  From among this assortment, the customer chooses a product according to some known, general model of choice.  If the customer chooses a product, the system earns a customer-type-dependent and product-dependent reward, and the inventory of the corresponding resource is depleted by one.  The system's goal is to maximize the total expected reward that it earns over the horizon.

The problem above has application in a number of settings.  In the revenue management of parallel flights, the resources are flight legs that share the same origin and destination.   Each product is a ticket, which is defined by a flight leg, an associated fare, and a set of purchase restrictions or ancillary services for the corresponding fare class.  A customer arriving at the system chooses dynamically among the assortment of tickets offered.  The type of a customer can be based on information such as the customers' past purchase history, their arrival time,  the parameters of the customers' search query, whether or not they belong to a loyalty program, their level within the loyalty program,  their location, and any other information that can be observed by the system at the time of purchase.

In assortment-planning problems, the products and resources are one and the same.  In this case, the type of a customer can be based on information similar to above, and whether they belong to a special program such as Amazon's Prime program, Bloomingdale's Loyallist rewards program, or Sephora's Rewards Boutique program.  Rewards are adapted to the customer type to capture a combination of personalized prices and the value of serving the customer type in the long term.  For example, Target is known to send coupons to shoppers that it identifies as potential expecting mothers, because data shows that important life events, such as birth, can change people's shopping habits.  Target finds it beneficial in the long run to favor these shoppers in order to induce the habit to buy at Target \citep{duhigg2012companies}.

In web- and mobile-based self-scheduling systems such as ZocDoc, the resources and products are the same and correspond to appointment slots that take place at different times, at various clinic locations, with various physicians. The patient type can be defined with reference to a patient's status as new or existing patient, their insurance, the nature of concern, etc.
The reward of assigning a slot to a patient can be chosen to capture a combination of the revenue expected for the visit and the patient's priority.

%
%

The above class of problems are closest to choice-based revenue-management problems that have attracted intense interest over the past ten years. Researchers have focused on stochastic dynamic formulations of these problems.
Gallego, Iyengar, Phillips, and Dubey (2004)\nocite{gallego2004managing} are among the first to introduce a consumer choice model to the network revenue-management literature. The choice based linear program (CDLP)  proposed by them has been widely used to approximate the stochastic dynamic optimal solution. They show that the optimal value of the stochastic dynamic problem is bounded above by the optimal solution to the CDLP, and further, that it approaches this upper bound asymptotically.  There are various ways of using the CDLP solution to obtain practical policies, such as the bid price heuristic of Talluri and van Ryzin (2005) \nocite{talluri2005theory} and the dynamic programming decomposition approach of Liu and van Ryzin (2008)\nocite{LiuV2008}.  Subsequently tighter approximation methods have been proposed to improve upon the CDLP upper bound, \citep{zhang2009af,  Meissner2012459, KunnumkalT2010}.  However, little is known about the theoretical performance of these methods outside of the asymptotic regime, and the tighter bounds usually come with significant computational cost.

In this paper, we introduce the first algorithms with theoretical performance guarantees for the above class of problems.  Our approach offers both modeling flexibility, ease of implementation, and theoretical performance characterization.  Several features of our model merit attention:
\begin{description}
\item [Personalization.]  We allow multiple customer types to be modeled. The customer types are based on information that can be observed by the system at the time of purchase. We offer assortments of products that are customized to the customer types.  Personalization has been widely used by companies that collect data on customer characteristics. It has been used, for example, by Amazon to recommend products to customers based on their purchase history; by Groupon, Yelp, and Foursquare to offer discounts to customers based on their location; and by grocery stores to offer customers coupons based again on their purchase history (Golrezaei, Nazerzadeh and Rusmevichientong 2014\nocite{golrezaei2014real}).  As we shall discuss, the literature on personalized revenue management is still limited.


\item [Substitution across fare classes.]  We allow for simultaneous consideration of multiple products that might correspond to multiple fare classes, with no restriction on substitution behavior.  Thus, we are able to capture the substitution of capacity across fare classes in revenue-management applications. Although the substitution effect has been studied, it has been limited to within fare classes  (Zhang and Cooper 2005\nocite{zhang2005revenue}). This assumption turns out to be rather restrictive.  
Dai, Ding, Kleywegt, Wang and Zhang (2014) \nocite{JimDai2014} have found in recent analysis of empirical airline data that customers have much greater demand for the cheapest alternative than for the second cheapest alternative even when the price difference is small.

\item [Substitution across time.]  We model multiple resources with expiry dates that can fall within the horizon.  Thus, we are able to explicitly model inter-temporal demand substitution.  Indeed, our model is among the first in the revenue-management literature to \emph{explicitly} capture inter-temporal substitution.  Previous models have usually assumed that the expiry dates fall beyond the horizon.  Thus, they can only capture intertemporal substitution in an implicit manner, by changing the demand arrival process or the selection probabilities over time.

\item [Non-stationarity.]  We allow demand arrivals to be non-stationary and stochastic over the horizon.  Our performance guarantees are robust to dramatic changes in the demand rate. Past approximations, such as that of Liu and van Ryzin (2008)\nocite{LiuV2008}, can be extended to the environment with time varying arrivals if the demand varies slowly over time, but perform poorly when the demand is volatile. Many of the remaining prior methods \citep{zhang2009af,  Meissner2012459, KunnumkalT2010} have assumed stationary arrival. 

\item [Dynamic product substitution.] We dispense with the commonly-used static-substitution assumption, which implies that a customer, who finds that his request is stocked out, will leave the system forever. We dynamically adjust our offered assortment according to current inventory levels, current demand type and future expected demands.  We include only products with positive inventory levels in our assortments. For this reason, the customer behaviors that our model captures are much more realistic.
\end{description}

Our contributions are as followed.
\begin{itemize}

\item We propose the first column-generation algorithms to solve the CDLP in settings where this problem is NP-hard.  Our algorithms can generate $\epsilon$-optimal solutions to the CDLP for any given $\epsilon>0$.  It is applicable to a variety of choice models, including the Mixed Multinomial Logit (MMNL), which can approximate any random utility model to any precision level \citep{mcfadden2000}. Our algorithms build on existing polynomial-time approximation schemes (PTAS) or fully polynomial-time approximation schemes (FPTAS) that can approximately solve the underlying assortment-planning sub-problems.

\item  We derive theoretical performance characterization for our general class of choice-based resource-allocation problems.  We prove that our algorithms are guaranteed to produce an expected reward no worse than $\frac{1}{2}(1-\epsilon)$ times that of $OFF$, where $\epsilon$ is the error in computing an optimal solution to the CDLP, and $OFF$ is an optimal offline algorithm that knows a priori the realization of all demand arrivals and makes optimal decisions given this information.  To the best of our knowledge, this is the first constant relative performance characterization for this class of problems. Our algorithms are highly intuitive and simple to implement.

\item We prove that $\frac{1}{2}$ is the best possible constant ratio that can be achieved between the expected reward of an online algorithm and the expected reward of  $OFF$.  Thus, our algorithms achieve within $(1-\epsilon)$ of the best possible constant relative performance.  They achieve the upper bound of $\frac{1}{2}$ for a variety of choice models for which the CDLP can be solved exactly.
  
\item We show that our algorithm has expected total reward at least as high as that of the deterministic algorithm proposed by Gallego, Iyengar, Phillips, and Dubey (2004) \nocite{gallego2004managing}.  We also prove that this classical algorithm has expected reward that is no worse than $1/e=0.368$ that of $OFF$.  Thus, we prove that our algorithms are asymptotically optimal, therefore competitive with existing algorithms, according to the main theoretical performance characterization known prior to our work.

\end{itemize}

\section{Literature Review}

We will summarize the streams of works that are most closely related to our paper, many of which are in revenue-management.  We refer the reader to Talluri and van Ryzin (2005) \nocite{talluri2005theory}for a comprehensive review of the larger revenue-management literature. 

\subsection{Single-Leg Revenue Management}
Inspired by the news-vendor problem, Littlewood (1972) \nocite{littlewood2005special} describes a simple technique for setting a booking limit for low-fare tickets when there are two fare classes and a single-leg flight. Later, Belobaba (1989)\nocite{belobaba1989} extends the model to the case with multiple booking classes and proposes an expected marginal seat revenue (EMSR) heuristic. Subsequently, a framework for determining booking limits for a single-leg flight with mutually independent demand classes that arrive in sequential blocks is developed. The earliest dynamic model is perhaps attributable to Mayer (1976)\nocite{mayer1976}, who introduces a dynamic-programming model and compares it with the static control derived from Littlewood's rule. Dynamic model allows for a more precise formulation of the customer arrival process. For example, the low-to-high-fare arrival assumption can be relaxed. Lee and Hersh (1993)\nocite{lee1993model} consider a general multiclass dynamic seat-allocation model with non-stationary demand. They show that the optimal policy is a monotone-threshold policy.  That is, a ticket is open for purchase if and only if its fare is no smaller than the expected marginal revenue over the remaining booking horizon. Talluri and van Ryzin (2004)\nocite{TalluriV2004} analyze a a single-leg revenue management problem with customer choice. They show that an optimal policy can be characterized by an ordered sequence of ``efficient" offer sets, which can provide the most favorable trade-off between expected revenue and expected capacity consumption. They give conditions under which the efficient sets have a nested structure.

Within this literature, a number of papers address the design of policies for revenue management that are robust to the distribution of arrivals. \citet{ball2009toward} analyze online algorithms for the single-leg revenue-management problem. Their performance metric is the traditional competitive ratio that compares online algorithms with optimal offline algorithms under the worst instance of demand arrivals. They prove that the competitive ratio cannot be bounded by any constant when there are arbitrarily many customer types. In our work, we relax the definition of competitive ratio, and show that our algorithms achieve a constant competitive ratio (under our definition) for any number of customer types and for a more general multi-resource model.  Qin, Zhang, Hua, and Shi (2015)\nocite{CongApproximationRM} study approximation algorithms for an admission control problem for a single resource when customer arrival processes can be correlated over time.  They use as the performance metric the ratio between the expected cost of their algorithm and that of an optimal stochastic dynamic algorithm. Our performance metric is stronger than theirs as we compare our algorithms against an optimal offline algorithm, instead of an optimal stochastic dynamic algorithm.  Qin, Zhang, Hua, and Shi (2015)\nocite{CongApproximationRM} prove a constant approximation ratio for the case of two customer types, and also for the case of multiple customer types with specific restrictions. In addition, they allow only one type of resource to be allocated.  In our model, we assume arrivals are independent over time, but we allow for multiple customer types and multiple resources without additional assumptions.

\subsection{Online Matching}
Our work is closely related to works on online bipartite matching problems.  In these problems, the set of available resources is known and corresponds to one set of nodes.  Demand requests arrive one by one, and correspond to a second set of nodes.  As each demand node arises, its adjacency to the resource nodes is revealed.  Each edge has an associated weight.  The system must match each demand node irrevocably to an adjacent resource node.  The goal is to maximize the total weighted or unweighted size of the matching.  There is no choice process that is modeled.

The online unweighted bipartite matching problem is originally shown by Karp, Vazirani and Vazirani (1990)\nocite{karp1990optimal} to have a best competitive ratio of $0.5$ for deterministic algorithms and $1-1/e$ for randomized algorithms. Our work generalizes the online weighted bipartite matching problem.  When demands are chosen by an adversary, the worst-case competitive ratio of this problem cannot be bounded by any constant \citep{mehta2012online}. Many subsequent works have tried to design algorithms with bounded performance ratios for this problem for more regulated demand processes.

Specifically, three types of demand processes have been studied. The first type of demand processes studied is one in which each demand node is independently and identically chosen with replacement from a \emph{known} set of nodes. Under this assumption, \citet{jaillet2013online, manshadi2012online, bahmani2010improved, feldman2009online} propose online algorithms with competitive ratios higher than $1-1/e$ for the unweighted problem. Haeupler, Mirrokni, Vahab and Zadimoghaddam (2011)\nocite{haeupler2011online} study online algorithms with competitive ratios higher than $1-1/e$ for the weighted bipartite matching problem.

The second type of demand processes studied is one in which the demand nodes are drawn randomly without replacement from an unknown set of nodes.  This assumption has been used in the secretary problem (Kleinberg 2005, Babaioff, Immorlica, Kempe, and Kleinberg 2008)\nocite{kleinberg2005multiple,babaioff2008online}, ad-words problem \citep{goel2008online} and bipartite matching problem (Mahdian and Yan 2011, Karande, Mehta, and Tripathi 2011)\nocite{mahdian2011online,karande2011online}.

The third type of demand processes studied is one in which each demand node requests a very small amount of resource. This assumption, called the \emph{small bid} assumption, together with the assumption of randomly drawn demands, lead to polynomial-time approximation schemes (PTAS) for problems such as ad-words \citep{Devanur09theadwords}, stochastic packing (Feldman, Henzinger, Korula, Mirrokni, and Stein 2010)\nocite{feldman2010online}, online linear programming (Agrawal, Wang, Zizhuo and Ye 2009)\nocite{agrawal2009dynamic}, and packing problems \citep{molinaro2013geometry}. Typically, the PTAS proposed in these works use dual prices to make allocation decisions. Devanur, Jain, Sivan, and Wilkens (2011)\nocite{devanur2011near} study a resource-allocation problem in which the distribution of nodes is allowed to change over time, but still needs to follow a requirement that the distribution at any moment induce a small enough offline objective value. They then study the asymptotic performance of their algorithm. In our model, the amount capacity requested by each customer is not necessary small relative to the total amount of capacity available.  Therefore, the analysis in these previous works does not apply to our problem.

Our work builds upon recent results by Wang, Truong and Bank (2015)\nocite{wangTB2015}, who propose online algorithms with competitive ratio of $0.5$ for the bipartite matching problem with heterogeneous demands.  They assume that demands arrive according to non-homogenous Poisson processes.  For this class of problems, they show that the bound of $0.5$ achieved by their algorithm is tight.  Our model extends theirs to allow for endogenous customer behavior, and decisions that control assortments of products to offer, rather than those that match customers directly to resources.

\subsection{Choice-based network revenue management}
Our work extends the literature on choice-based network revenue management (CNRM).  These models assume that there is some fixed, finite amount of resources, for example, flight legs.  Products are sold, which are built from one or more resources.   The models focus on dynamically adjusting the set of offered products as a function of remaining capacity of the resources and remaining time in the selling horizon.  They assume that consumer demands are dependent on the set of offered products. 

Our algorithms build on the use the CDLP that was first introduced by Gallego, Iyengar, Phillips, and Dubey (2004)\nocite{gallego2004managing}. They propose an efficient column generation algorithm to solve the CDLP and suggest a deterministic control policy that  uses the primal optimal solution as the proportion of time to offer each assortment.  More recently, Gallego, Ratliff and Shebalov (2014) \nocite{gallego2014general}reformulate the CDLP as a polynomial-size sales-based LP under a general class of attraction models. Liu and van Ryzin (2008)\nocite{LiuV2008} give a dynamic-programming formulation for CNRM. They come up with a heuristic algorithm by decomposing over flight legs, with the opportunity cost of each flight leg being generated from the dual values of the capacity constraints in the LP. Kunnumkal and Topaloglu (2010)\nocite{KunnumkalT2010} propose another dynamic-programming-decomposition algorithm for CNRM; they allocate revenue associated with itinerary among the different flight legs and solve a single-leg revenue-management problem for each leg. Zhang and Cooper (2005) \nocite{zhang2005revenue} dynamically control the inventory of parallel flights for a common itinerary while assuming demands are substitutable across different schedules, but not across different fare classes. Zhang and Adelman (2009)\nocite{zhang2009af} present an affine approximation to the value function. Jasin and Kumar (2012)\nocite{JasinKumar2012} study the performance of the heuristic of resolving the deterministic LP periodically while setting all random variables at their expected future values.  They provide both upper bound and lower bound for the expected revenue loss compared to the optimal policy.  All of the above papers assume a single market segment, or multiple market segments with disjoint consideration sets. They also assume that the customer type is not observable to the seller, so the seller must decide on a common assortment for all segments.

As Bront, Mendez-Diaz, and Vulcano (2009)\nocite{BrontMV2009} point out, when customers belong to overlapping segments, even solving for the CDLP is NP-hard.  They propose a heuristic column-generation algorithm for solving the CDLP. Our model is different in that we can observe the customer type and can provide personalized assortments.

Jaillet and Lu (2012)\nocite{jaillet2011online} design near-optimal learning-based online algorithms for dynamic resource-allocation problem. They do not model consumer choice.  Further they assume that demand arrivals are stationary, and that the amount of resource used by each demand unit is very small relative to the capacity.  In contrast, we will model consumer choice, allow non-stationary demands, and allow individual demand requirements to be large relative to capacity.

\subsection{Revenue management of parallel flights}
Zhang and Cooper (2005)\nocite{zhang2005revenue} consider the seat-inventory control of multiple parallel single-leg flights in the presence of dynamic customer-choice behavior. Each customer's choice is modeled by his preference mapping, and they assume the preferences only shift among different flights, not across fare classes. They use lower and upper bounds to approximate the expected revenue from the optimal stochastic programming, and suggest control policies based on the approximated marginal value. Dai, Ding, Kleywegt, Wang and Zhang (2014) \nocite{JimDai2014} describe a revenue-management problem of a major airline that operates in a very competitive market involving two hubs and having more than 30 parallel daily flights. They observe demand discontinuities from the industrial date, i.e, demand spikes for the cheapest available fare classes and for fully refundable fare classes. They incorporate the discontinuity into a logit model and focus on a deterministic formulation of the problem. In a departure from previous literature, they take the competitor's response into consideration by modeling the attractiveness of the no-purchase option with a random coefficient. They also show that under some conditions, the CDLP can be solved efficiently by column generation.

\subsection{Assortment planning}
Assortment planning problems are special class of revenue-management problems in which the products sold are the resources themselves.  That is, the products are not built from one or more simpler resources.  In assortment planning, a retailer needs to decide the set of products to offer at various times over a selling horizon.  Usually this decision is jointly considered with inventory decisions.

Earlier works considered the static assortment problem.  In these problems, customers have no knowledge about the status of inventory.  They make purchase decision only based on the offered product set.  If their selection is stocked out, a second choice will not be made.  Van Ryzin and Mahajan (1999)\nocite{VanRyzinM1999} show that an optimal assortment is composed of the most popular products when all the products have the identical price and cost and customers' demands are governed by the Multinomial Logit model of choice. Subsequent literature has considered various choice models and prices and costs structures. Cachon, Terwiesch and Xu (2005)\nocite{P.Cachon2005} use a more general choice model in which they incorporate search costs and shows that ignoring customer search will lead to less assortment variety since in equilibrium the seller needs a bigger sized assortment to attract more customers. Topaloglu (2013)\nocite{Topaloglu2013stock} works on a similar problem as van Ryzin and Mahajan (1999)\nocite{VanRyzinM1999}, but instead of deciding a single offering set, he allows multiple assortments and decides the duration of time that each is offered.

Our paper features a dynamic, or stockout-based model of consumer choice.  In these models, customers base their choice from the products that are in stock upon their arrival. Mahajan and van Ryzin (2001)\nocite{MahajanV2001} choose initial inventory levels to maximize the expected profit under dynamic substitution.  However, they show the objective function is not even quasiconcave. Therefore, they suggest a stochastic-gradient algorithm for the problem. Kok and Fisher (2007)\nocite{kokfisher2007} propose an estimation method to obtain both the original demand rate and the substitution rate.  They also present an iterative heuristic to solve the problem of joint assortment planning and inventory optimization. Honhon, Gaur and Seshadri (2010)\nocite{honhonor2010} consider multiple customer types.  They derive some structural properties for the optimal assortment.  

Within the literature on dynamic consumer choice, our paper is closest to the literature on personalized dynamic assortment planning.  These assortments are dynamically optimized depending on inventory levels, and tailored to the customer segment, provided that the system can observe the customer's segment and knows the preferences of each segment. Bernstein, Kok and Xie (2010)\nocite{Bernstein10dynamicassortment} are the first to propose the idea of assortment customization in presence of heterogeneous customer segments. To obtain some structural properties of the optimal policy, they assume all products are not functionally differentiated and have the same price; and they propose a heuristic by implementing a newsvendor-type approximation to the marginal revenue of each product. \cite{chan2009stochastic} relax the restrictive uniform price constraint and show that a myopic policy achieves at least 0.5 the expected revenue of the stochastic optimal policy in non-stationary settings.  Golrezaei, Nazerzadeh and Rusmevichientong (2014)\nocite{golrezaei2014real} extend this result by showing that under adversarially chosen demand there is an online algorithm that achieves at least 0.5 the cost of an optimal offline policy in the worst case.

Similar to Golrezaei, Nazerzadeh and Rusmevichientong (2014)\nocite{golrezaei2014real}, we propose simple and robust online algorithms for dynamically determining the set of offered products, based on the inventory of products available and the time that remains in the horizon.  However, our model captures rewards that depend on the customer type and not just on the products sold.  In revenue-management applications, these rewards capture differentiated fares for the same seat capacity.  In assortment planning applications, our rewards capture personalized prices or discounts, and differentiated rewards for serving different customer groups, such as Amazon Prime customers versus regular customers.  In self-scheduling systems, our rewards capture differentiated priorities among different customer groups such as urgent versus non-urgent patients, new versus existing patients, and regular versus follow-up visits, etc.  In the context of online matching, the models of Golrezaei, Nazerzadeh and Rusmevichientong (2014)\nocite{golrezaei2014real}, \cite{chan2009stochastic}, and others in the assortment-planning literature extend the \emph{vertex-weighted bipartite matching problem}, where the revenue earned is a function of only the resource nodes that are used.  This problem is a special case of the \emph{edge-weighted bipartite matching problem} that we extend.  In these more general problem, the revenue earned is a function of both the resource nodes that are used, and the demand nodes that they are matched with.  

\subsection{Appointment Scheduling with Choice}
Our work is related to the literature on appointment scheduling \citep{guerriero2011operational,may2010surgical,cardoen2010operating,gupta2007surgical}.  Patient preferences are an important consideration in many real scheduling systems. In the literature considering patient preferences, \citet{gupta2008revenue} consider a single-day scheduling model where each arriving patient picks a single slot with a particular physician.  The clinic accepts or rejects the request. Our model generalize their framework to a multi-period setting. We also characterize the theoretical performance of algorithms in an online setting, whereas they use stochastic dynamic programming as the modeling framework and develop heuristics.

A multi-day, single-patient-type, stationary model has been studied by Feldman, Liu, Topaloglu and Ziya (2014)\nocite{feldman2014appointment}.  They assume that patients have preferences for slots that can be captured by the multinomial logit model.  This model is essentially a dynamic assortment-planning model.  They derive structural results for the optimal policy, and exhibit a heuristic that is asymptotically optimal. In contrast, we characterize the theoretical performance of our algorithms in both asymptotic and non-asymptotic regimes.  We also model heterogenous patients and non-stationary demand, both of which features are especially important in healthcare settings, where patients frequently have differing priorities and demands are usually non-stationary (Wang, Truong and Bank 2015)\nocite{wangTB2015}.  Finally, we allow a very general model of patient choice to be used.

Our model closely follows the previously discussed model of Wang, Truong and Bank (2015)\nocite{wangTB2015} for multi-day, multi-patient-type settings.   Their model is useful in applications in which the system can observe \emph{revealed} patient preferences before assigning appointment slots to them.  Requiring patients to explicitly specify their preferences  can be a cumbersome exercise.  Therefore, some systems such as ZocDoc make scheduling more user-friendly by offering sets of open slots to patients and allowing them to choose.  Our model captures the latter setting.  Thus, our contribution is to add the control of assortments of slots to offer to the model of Wang, Truong and Bank (2015)\nocite{wangTB2015}.

\section{Model}

We consider a continuous horizon $[0,1]$.  There are $L$ different resources, $N$ different products, and $K$ customer types.  Each resource $l$, $l=1,\ldots,L$, has a capacity $C_l$ and an expiration time $t_l$.  Each product $n$, $n=1,\ldots,N$, consists of a single resource $l_n$.  A product might be a resource offered with a specific price and a certain set of restrictions.  For expositional simplicity, each product $n$ earns a reward of $r_n$ regardless of the customer type served.  However, the reward can be made to depend on the customer type served without changing any of the results that follow. In the more general case, the reward would be indexed by both the product and the customer type.

Customers of type $k$ arrive according to a non-homogenous Poisson process with {\it known} rate $\lambda_k(t)$. We assume that the customer type is observable by the system.  When a customer of type $k$ is presented with an assortment $S$, the customer will choose a product from the assortment following a {\it general} choice model. Let $P^k(n,S)$ be the probability that a customer of type $k$ chooses product $n$ from assortment $S$. Note that \emph{each customer type can be modeled with a different choice model.} We assume that the no-purchase option, or product $0$, is included in each assortment, with $r_0=0$.


When a customer of type $k$ arrives, the system must decide which assortment to offer to the customer. If the customer purchases product $n$ from the offered assortment, the reward $r_n$ is earned and the inventory of resource $l_n$ is depleted by $1$. The objective of the system is to maximize the expected total reward.

Let $c = (c_1,c_2,...,c_L)$ denote the vector of remaining inventory, where $c_l \in \{0,1,...,C_l\}$ represents the remaining inventory of resource $l$. Under a policy $\Pi$, let $V^\Pi(c,t)$ denote the expected future reward and $A_k^\Pi(c,t)$ denote the assortment offered to customers of type $k$ at time $t$ with inventory $c$. In general, $A_k^\Pi(c,t)$ can be random.

Let $N_l \equiv \{n \in \{1,2,...,N\}: l_n = l\}$ be the set of products that consist of resource $l$. Let $e_l$ be the unit vector with the $l$-th position being $1$. The dynamics of $V^\Pi(c,t)$ is governed by
\begin{equation}\label{eq:Dynamics} 
\frac{\partial V^\Pi(c,t)}{\partial t} = - \sum_{k=1}^K \sum_{l=1}^L  R_t^{kl}(A_k^\Pi(c,t), \Delta_l V^\Pi(c,t)),
\end{equation}
where
\begin{equation}\label{eq:DeltaV}
 \Delta_l V^\Pi(c,t) \equiv \left\{ \begin{array}{ll} V^\Pi(c,t) - V^\Pi(c-e_l,t) &\text{ if } c_l>0\\ \infty & \text{ if } c_l =0\end{array} \right.
 \end{equation}
is the marginal value of expected future reward with respect to resource $l$, and
\begin{equation}\label{eq:RewardRate}
 R_t^{kl}(S,z) \equiv \lambda_k(t) \sum_{n \in N_l} \mathbf{E}[\mathbf{1}(n \in S)  P^k(n,S)] (r_n - z)
 \end{equation}
is the rate at which the expected future reward of resource $l$ changes due to customers of type $k$. Here the expectation in (\ref{eq:RewardRate}) is taken over $S$. $V^\Pi(c,t)$ must satisfy the boundary conditions $V^\Pi(0,t) = 0$ and $V^\Pi(\cdot,1) = 0$.

The total expected reward of policy $\Pi$ is $V^\Pi(C,0)$, where $C =(C_1,C_2,...,C_L)$ is the vector of initial inventory.

\subsection{Action space} 
Existing models differ in whether an assortment is allowed to contain products with zero inventory. A model assumes \textit{static substitution} if an assortment can contain any product.  A customer who chooses a product with zero inventory leaves the system without affecting the total reward. A model assumes \textit{dynamic substitution} if assortments must not contain products with zero inventory.

Dynamic substitution is much more realistic than static substitution, although dynamic substitution requires more complex analysis. Recent literature has mostly focused on dynamic substitution.

Our paper assumes dynamic substitution. The formal definition is as follows.
\begin{assumption}
For any given state $(c,t)$, if $c_l =0$, we must have $N_l \cap A_k^\Pi(c,t)=\emptyset$.
\end{assumption}

We will propose a policy that achieves the best performance guarantee under this assumption. To motivate the idea of the policy, however, we will first analyze some intermediate policies that generate performance bounds under the static-substitution assumption. We will show that our final policy that works under dynamic substitution dominates these intermediate policies.

\subsection{More on time-dependent  effects}
Note that in our model, the resources might perish over the horizon.  That is, it is possible for $t_l < 1$ for certain resources $l$.  We will define customer types to ensure that any customer $k$ arriving after time $t_l$ has selection probability $P^k(n,S)=0$ for any product $n$ that is based on $l$, and any assortment $S$ that includes $n$.  As time moves forward, the types of customers who are arriving with positive probabilities will change in our model in order to capture any time-dependent changes in demand behavior.

The horizon is also finite in our model.  As the end of horizon approaches, there might be products that expire after the end of horizon entering into customers' consideration sets.  In the language of revenue management, as the end of horizon approaches,  customers might be increasingly choosing to purchase flights that depart after the end of the horizon.  Existing models can address this effect by changing the attractiveness of the no-purchase option as the end of horizon nears.  A similar strategy can be applied here.  However, our model can be used to capture inter-temporal substitution more explicitly near the end of the horizon as follows.  For concreteness, we explain the strategy in the language of revenue management.
\begin{itemize}
\item Include consideration sets with flights that depart up to time $T + \Delta$, where $\Delta$ might be a week. 

\item Allocate capacity $c_l$ to flights $l$ with departure time $t_l > T$ for sale during $[0,T]$. In this case $c_l$ may be less than the actual capacity for flight $l$. The quantity $c_l$ may be set by management or may be the solution to a higher-level optimization problem. As an example, if the capacity of a flight that departs between $T$ and $T+ \Delta$ is 100, we might allocate 85 seats to the optimization problem over $[0,T]$ if we estimate that we can sell 15 seats at a higher price during the interval $[T, T+\Delta]$. 
\end{itemize}
By including post-horizon flights and some of their capacity, we will be able to incorporate most of the inter-temporal substitution behavior explicitely into our model.

\subsection{Competitive ratio}
It is practically impossible to compute the optimal policy for our problem because of the ``curse of dimensionality'' in both the state and decision spaces. Instead, our goal is to give an online algorithm with expected total reward that is bounded by a constant factor of an optimal offline policy $OFF$. 

A policy $\Pi$ is online if the decision $A_k^\Pi(c,t)$ is adapted to the information up to time $t$, including the future arrival rates $\lambda_k(t)$, for $k=1,2,...,K$ and $t \in [0,1]$, that are known a priori. On the other hand, the optimal offline decision depends on the information of all realizations of future arrivals over the horizon, but the randomness in customer choices is still exogenous to \textit{OFF}. We will also require that $OFF$ follows the dynamic-substitution assumption in the sense that it never offers a product with zero inventory. 

Let $V^{OFF}$ denote the expected reward of $OFF$. Note that $V^{OFF}$ does not need to satisfy the dynamic equation (\ref{eq:Dynamics}) for online algorithms. 
The following is our definition of competitive ratio.
\begin{definition}
An algorithm $\Pi$ is \emph{$\alpha$-competitive} if its total expected reward $V^\Pi(C,0)$ satisfies
\[ V^\Pi(C,0) \geq \alpha \cdot V^{OFF} \]
under any values of $r_{n}, \nu^k_n, \lambda_k(t), C_l$, for $n=1,2,...,N;\ k = 1,2,...,K;\ l = 1,2,...,L;\ t\in [0,1]$. 
\end{definition}

\section{Choice-Based Deterministic Linear Program}
Before introducing our algorithms, we first characterize an upper bound on $V^{OFF}$ using a widely-used choice based deterministic linear programming (CDLP) 
formulation proposed by \citet{gallego2004managing}. \citet{LiuV2008} have shown that the CDLP is an upper bound on the expected reward of an optimal stochastic policy. In this paper, we prove a stronger result that the CDLP is an upper bound on the expected reward of an optimal offline policy \textit{OFF}.  Our algorithms and bounds will build on this CDLP.

Let $\cal{S}$ denote the set of all assortments.  Let $\Lambda_k \equiv \int_0^1 \lambda_k(s) ds$ denote the average total number of arrivals of type $k$ customers over the horizon. In the following choice-based CDLP, the decision $x_k(S)$ represents the probability of showing assortment $S \in \cal{S}$ to a type-$k$ customer upon his arrival.  Note that this quantity is independent of the time of that arrival.

\begin{align}
\begin{split} \label{eq:ChoiceBasedLP}
V^{CDLP}=& \max \sum_{k=1}^K \sum_{S\in \cal{S}} \Lambda_k x_k(S) \sum_{n\in S} P^k(n,S) r_n \\
\text{s.t. } &  \sum_{k =1}^K   \Lambda_k \sum_{n \in N_j} \sum_{S \in \cal{S}: n \in S} x_k(S) P^k(n,S) \leq C_j , \,\,\, \forall j = 1,2,...,L,\\
 & \sum_{S \in \cal{S}} x_k(S) \leq 1, \,\,\, \forall k =1,2,...,K;\\
 & x_k(S) \geq 0, \,\,\, \forall k =1,2,...,K,\ \forall S \in \cal{S}.
 \end{split}
\end{align}

\begin{theorem} \label{thm:upperbound}
$V^{CDLP}$ is an upper bound on $V^{OFF}$.\end{theorem} 

Our algorithms rely on an optimal solution of the CDLP.
 Since there are $L+K$ constraints, we know that at optimality, the solution involves at most $L+K$ different assortments with positive displaying probabilities. However, with $K2^N$ variables, this CDLP could be very difficult to solve in practice. As suggested by Gallego, Iyengar, Phillips, and Dubey (2004)\nocite{gallego2004managing}, column generation techniques can be used to circumvent these difficulties. 
 
\subsection{Solving the CDLP by Column Generation}
Our algorithm requires a solution to the CDLP.  We propose solving the CDLP by column generation.  Column generation works by expanding a set $\cal{H}^d$ of \emph{active} assortments, i.e., those with positive displaying probabilities, in each iteration $d$. We start with a limited number of columns, namely $\cal{H}^{1}$.  We then solve a reduced LP that involves only these columns. Conditional on the dual values from the reduced LP in the current iteration, we then calculate the reduced cost of potential assortments that have not yet been considered. If the reduced cost is strictly positive, we incorporate the corresponding column into a new active set and re-optimize. If there is no potential assortment with positive reduced cost, the current solution is optimal. A detailed illustration is as follows:
\begin{enumerate}
\item Initially, let $\cal{H}^1$ be any collection of assortments. Set $d\leftarrow 1$.
\item Solve the following reduced linear program, in which $x_k(S)$ is the decision variable, for all $k = 1,2,...,K$ and for all $S \in \cal{H}^{d}$.
\begin{align*}
\begin{split}
 V^d \equiv \max &\,\,\,\,\,\,\, \sum_{k=1}^K \sum_{S \in \cal{H}^{d}} \Lambda_k x_k(S) \sum_{n \in S} r_n P^k(n,S)\\
\text{ s.t. } & \sum_{k=1}^K \Lambda_k \sum_{n \in N_j} \sum_{S\in \cal{H}^{d}: n \in S}  x_k(S) P^k(n,S) \leq C_j, \,\,\forall j = 1,2,...,L\\
& \sum_{S \in \cal{H}^{d}} x_k(S) \leq 1, \,\,\, \forall k = 1,2,...,K\\
& x_k(S) \geq 0, \,\,\, \forall k = 1,2,...,K,\ \forall S \in \cal{H}^{d}.
\end{split}
\end{align*}
Let $\pi(1), \pi(2), ..., \pi(L), \sigma(1), \sigma(2),...,\sigma(K)$ be the optimal dual variables for the above reduced LP. The reduced cost corresponding to $x_k(S)$, which is the offering probability of assortment $S \in \cal{S}$ to consumer type $k$, is 
\[ \Lambda_k \sum_{n \in S, n>0} [r_n - \pi(l_n)] P^k(n,S)  - \sigma(k).\]
We must check whether there is any assortment $S \notin \cal{H}^{d}$ that has a strictly positive reduced cost. This can be done by solving the following optimization problem for each $k$ 
\begin{equation}\label{eq:CGSubproblem}
L_k\equiv \max_S \{ \Lambda_k \sum_{n \in S, n>0} [r_n - \pi(l_n)] P^k(n,S)  \}
\end{equation}
and compare the value with $\sigma(k)$. Note the above is just an assortment optimization problem with $r_n - \pi(l_n)$ being the profit of product $n$.

\item Let $k^* \in {argmax}_k L_k-\sigma(k)$, if $L_{k^*} -\sigma(k^*)\leq 0$, the solution of the current reduced LP is optimal for the primal CDLP. Stop.
\item Else, set $\cal{H}^{d+1} \leftarrow \cal{H}^{d} \cup \{S^d_{k^*}\}$, where $S^d_{k^*} \in {argmax}_S  \Lambda_{k^*} \sum_{n \in S, n>0} [r_n - \pi(l_n)] P^{k^*}(n,S)$ and set $d \leftarrow d + 1$
\end{enumerate}

Column generation is successful only if \eqref{eq:CGSubproblem} can be solved efficiently for all values of $(r,\pi)$. Clearly, this is impossible for arbitrary assignments of the choice probabilities $P(n,S), S \in \cal{S}$. 

\subsection{Solving the Column-Generation Subproblem}
Now let us analyze the complexity of each column-generation step. Gallego, Iyengar, Phillips, and Dubey (2004)\nocite{gallego2004managing} show that if the customer choices follow the structure
\[ P^k(n,S) = \frac{\mu_{kn} + \nu_{kn}}{\sum_{i=1}^N \mu_{ki} + \sum_{i\in S} \nu_{ki} + 1}, \,\,\,\,\,\forall n \in S,\]
then the column-generation subproblem can be solved exactly by a simple sort. The above class of models cover a wide range of choice models, for example, the Independent-Demands, Multinomial Logit (MNL), and General Attraction models (GAM). Recently, Gallego, Ratliff and Shebalov (2014) \nocite{gallego2014general} show that under the above class of choice models, the CDLP can be reformulated as a more compact sales-based linear program with only a polynomial number of variables.

There are other choice models under which the optimal assortment can be found in polynomial time, for example, the $d-$Level Nested Logit (Li, Rusmevichientong and Topaloglu 2015)\nocite{doi:10.1287/opre.2015.1355} and  Markov-Chain models (Blanchet, Gallego and Goyal 2013)\nocite{Blanchet:2013}. For those choice models, column generation can efficiently return the optimal value. However, under more complex choice models, the assortment problem might be NP-hard. For example, Bront, Mendez-Diaz, and Vulcano (2009)\nocite{BrontMV2009} show that the column generation sub-problem under the Mixed Multinomial Logit (MMNL) model is NP-hard. Moreover, when an optimal solution to the CDLP cannot be obtained, there is currently no algorithm that can generate an approximate solution to the CDLP that is guaranteed to be $\epsilon$-optimal for any given $\epsilon>0$, where $\epsilon$ is the error of the approximate solution defined as follows.

\begin{definition}
Let $x_k(S)$, for $k=1,2,...,K$ and $S \in \cal{S}$, be a feasible solution to (\ref{eq:ChoiceBasedLP}). We say that $x_k(S)$ is $\epsilon$-optimal if
\[ \sum_{k=1}^K \sum_{S \in \cal{S}} x_k(S) \sum_{n \in S} r_n P^k(n,S)\geq (1-\epsilon)V^{CDLP}.\]
\end{definition}

We propose the first column-generation algorithms specifically designed for the case when the sub-problem is NP-hard.  Our algorithm can generate $\epsilon$-optimal solutions to the CDLP, for any given $\epsilon>0$, for a variety of choice models, including the MMNL, which can approximate any random utility model to any precision level \citep{mcfadden2000}. Our algorithm builds on existing polynomial-time approximation schemes (PTAS) or fully polynomial-time approximation schemes (FPTAS) that can approximately solve the assortment-planning problem.

Compared with the column-generation algorithm defined before, we make a slight change in each iteration step. Instead of searching for the optimal assortment $S_k^*(\pi) \in \cal{S}$ for all $k$, which is NP-hard, we aim to find an assortment $S_k^d$ that satisfies
\[ \sum_{n \in S_k^d, n>0} P^k(n,S_k^d) [r_n - \pi(l_n)] \geq (1 - \frac{\epsilon}{1 + \epsilon}) \sum_{n \in S_k^*(\pi), n>0} P^k(n,S_k^*(\pi)) [r_n - \pi(l_n)] .\]
This modified problem can be solved in polynomial time using a PTAS or FPTAS that applies to a broad class of choice models, including the MMNL \citep{doi:10.1287/opre.1120.1093}.
The algorithm terminates when at the $\epsilon$ precision level, we cannot find any assortment $S_k^d$ such that
\[ \Lambda_k \sum_{n \in S_k^d, n>0} [r_n - \pi(l_n)] P^k(n,S_k^d)  - \sigma(k) >0 ,\]

\begin{theorem}\label{thm:nearOptimalCG}
For any given $\epsilon > 0$, the above column-generation algorithm generates an $\epsilon$-optimal solution.
\end{theorem}

\section{Upper Bound on the Competitive Ratio}

In this section, we show that $0.5$ is the best competitive ratio that any online algorithm can achieve for our model. This result implies later that our algorithms achieve the best performance guarantee.

\begin{theorem}
For the dynamic assortment-planning problem, no algorithm can achieve a competitive ratio higher than $0.5$.
\end{theorem}
\proof{Proof.}
Consider the following special case of our model. There is only one product and one associated resource. Then the problem becomes a single-resource reward management problem.  Wang, Truong and Bank (2015)\nocite{wangTB2015} have shown that the best competitive ratio for such single-resource problem is $0.5$ when the arrival rates are non-homogeneous. Therefore, for the general multi-product, multi-resource problem, the best competitive ratio is also $0.5$.
\halmos
\endproof

\section{The First-Come-First-Served Algorithm}
To motivate our main ideas, we first give an analysis of a simpler first-come-first-served ($FCFS$) algorithm. While $FCFS$ has a small proven performance guarantee under relaxed conditions, it leads to the full intuition behind our main algorithm.

For ease of exposition, we introduce the following somewhat impractical assumptions that relax the action space of online algorithms for the current section and Section \ref{sec:PrimalRouting}. Later in section \ref{sec:OPR}, we will show that these assumptions can be easily removed without loss of generality:
\begin{enumerate}
\item An assortment can contain products whose corresponding resource has zero inventory.  Customers' choices are unaffected by the inventory status of products.  If a customer chooses a product with zero inventory, this customer is treated as being rejected. Recall that this assumption is also called the \emph{static substitution} assumption. 
\item Even after a customer has chosen a product with positive inventory, an algorithm can still reject the customer, preventing the chosen product from being purchased.  That is, inventory will be unchanged after the interaction.  This assumption is used only in Section \ref{sec:PrimalRouting}. 
\end{enumerate}

The $FCFS$ algorithm is a naive implementation of an optimal solution to (\ref{eq:ChoiceBasedLP}):
\newcommand{\stepone}{(Pre-processing Step) Solve the CDLP (\ref{eq:ChoiceBasedLP}) to obtain an $\epsilon$-optimal solution $x^*_k(S)$ for $k = 1,2,...,K$ and $S \in \cal{S}$. Note that $\epsilon$ can be $0$ if it is possible to compute an exact optimal solution. For customers of type $k$, let $\cal{O}_k \equiv \{S \in \cal{S} : x_k^*(S) > 0\}$ be the set of assortments with positive displaying probabilities.}

\newcommand{\randomrouting}{(Random Offering Step) Upon an arrival of a type-$k$ customer at time $t$, offer an assortment $A_k^{Static}$ that is randomly picked from $\cal{O}_k$ such that $P(A_k^{Static} = S) = x_k^*(S)$ for all $S \in \cal{O}_k$. Note that $A_k^{Static}$ is independent of the state of inventory and time $t$.}  
\begin{enumerate}
\item \stepone

\item \randomrouting

\item ($FCFS$ Step) If the customer chooses a product $n \in A_k^{Static}$ and $n>0$, let the customer purchase it if the corresponding resource $l_n$ has positive remaining inventory. Otherwise, reject the customer.
\end{enumerate}

$FCFS$ is essentially a generalization of the deterministic algorithm of Gallego, Iyengar, Phillips, and Dubey (2004)\nocite{gallego2004managing} to a non-stationary, multi-customer-type environment.

Since $FCFS$ assumes static substitution, its expected future reward $V^{FCFS}(c,t)$ does not satisfy (\ref{eq:Dynamics}). Instead, the dynamic equation that $FCFS$ follows is
\begin{equation}\label{eq:DynamicsFCFS}
 \frac{\partial V^{FCFS}(c,t)}{\partial t} = - \sum_{k=1}^K \sum_{l=1}^L \mathbf{1}(c_l >0) R_t^{kl}(A_k^{Static}, \Delta_l V^{FCFS}(c,t)).
 \end{equation}
Note that the only difference between this equation and (\ref{eq:Dynamics}) is the additional term $\mathbf{1}(c_l >0)$, which implies that customers who choose products with zero inventory have no impact on the system.

Using an interesting analysis based on properties of Poisson processes, we can show that $FCFS$ already gives a constant performance guarantee under static substitution, as stated in the following theorem.
\begin{theorem} \label{thm:FCFS}
If $x^*$ is $\epsilon$-optimal then $V^{FCFS}(C,0) \geq \frac{1}{e}(1-\epsilon) V^{OFF}$.
\end{theorem}

\section{The Primal Routing Algorithm}
\label{sec:PrimalRouting}
The FCFS algorithm gives the following intuition. If we view the Random Offering Step as exogenous, the problem separates into $L$ single-resource revenue-management problems.  The demands arriving at each resource can be considered as {\it independent demands}.  

It is well-known that $FCFS$ is not optimal for the single-resource revenue-management problem, as it might not be optimal to always accept customers who arrive first. The following Primal Routing Algorithm (PR) optimally solves the single-resource revenue-management problem for each resource, by rejecting customers who want to purchase products at prices that are too low.  In this sense, $PR$ improves upon the performance of $FCFS$.  However, $PR$ has many disadvantages even compared to $FCFS$.  In particular, our analysis of $PR$ still assumes the \emph{static substitution} assumption that we made for $FCFS$.  In addition, $PR$ might offer a product as part of an assortment, then upon a customer's choosing the product, reject the customer.   This might happen even if the product has positive inventory.  

Nevertheless, the following performance analysis of $PR$ will help us to build intuition for the analysis of a more practical, but more intricate, algorithm that we will develop in Section \ref{sec:OPR}.

The difference between $PR$ and $FCFS$ lies in the Primal Routing Step that replaces the FCFS Step:

\newcommand{\primalrouting}{(Primal Routing Step) Let $c = (c_1,c_2,...,c_L)$ be the vector of inventory at time $t$. Suppose that a customer chooses a product $n \in A_k^{Static}$. Suppose that $n>0$ and $l = l_n$. Let the customer purchase the product if and only if $c_l>0$ and 
\[ r_n \geq \Delta_l V^{PR}(c,t) ,\]
where $V^{PR}(c,t)$ stands for the expected future reward of $PR$ given state $(c,t)$, and the definition of $\Delta_l V^{PR}(c,t)$ follows (\ref{eq:DeltaV}). The dynamic equation that $V^{PR}(c,t)$ must satisfy is
\begin{equation} \label{eq:DynamicsPR}
 \frac{\partial V^{PR}(c,t)}{\partial t} = - \sum_{k=1}^K \sum_{l=1}^L \cal R_t^{kl}(A_k^{Static}, \Delta_l V^{PR}(c,t)),
 \end{equation}
where
\[ \cal R_t^{kl}(S,z) \equiv \lambda_k(t) \sum_{n \in N_l} \mathbf{E}[\mathbf{1}(n \in S)  P^k(n,S)] [r_n - z]^+. \]
Note that $\cal R_t^{kl}(S,z)$ is different from $R_t^{kl}(S,z)$ defined in (\ref{eq:RewardRate}) as $\cal R_t^{kl}(S,z)$ applies the $[\cdot]^+$ operator to the difference $r_n-z$ between the reward of product $n$ and the cost value $z$.}


\begin{enumerate}
\item (Pre-Processing Step) Same as in FCFS.
\item (Random Offering Step) Same as in FCFS.
\item \primalrouting
\end{enumerate}

The dynamic equation (\ref{eq:DynamicsPR}) can be decomposed by resources as follows. For each resource $l$, consider the following single-resource reward function $V_l^{PR}(c,t)$ defined as
\begin{equation} \label{eq:dp}
 \frac{\partial V_l^{PR}(c,t)}{\partial t} = - \sum_{k=1}^K \cal R_t^{kl}(A_k^{Static}, \Delta_l V_l^{PR}(c,t))
 \end{equation}
 with boundary conditions $V_l^{PR}(0,t) = 0$ and $V_l^{PR}(c,1) = 0$.

It is easy to see that (\ref{eq:dp}) is just the Hamilton-Jacobi-Bellman equation that computes the optimal expected future reward for resource $l$ under arrivals over the horizon of $|N_l|$ demand classes, with demand class $n \in N_l$ bringing unit reward $r_n$ and arriving at the rate 
\[ \sum_{k=1}^K \lambda_k(t)  \mathbf{E}[\mathbf{1}(n \in A_k^{Static})  P^k(n,A_k^{Static})]. \]

This immediately leads to the following result
\begin{theorem}
The expected total reward of $PR$ is
\[ V^{PR}(C,0) = \sum_{l=1}^L V_l^{PR}(C,0).\]
In particular, the expected future reward that $PR$ earns from resource $l$, starting with state $(c,t)$, is $V_l^{PR}(c,t)$.
\end{theorem}
\proof{Proof.}
This theorem is a direct result of the properties of the Hamilton-Jacobi-Bellman equation (\ref{eq:dp}). \halmos
\endproof

Core to our competitive analysis of PR is the following property of the reward function.
\begin{theorem} \label{thm:SingleLegCR} For each resource $j=1,\ldots,L$,
\[ V_j^{PR}(C,0) \geq 0.5 \sum_{k=1}^K\Lambda_k \sum_{n \in N_j} \sum_{S \in \cal{O}_k : n \in S} x_k^*(S) P^k(n,S) r_n.   \]
\end{theorem}
The proof proceeds by proving that the expected revenue generated by each resource $j$ under $PR$ is at least half of what $OFF$ obtains.
Wang, Truong and Bank (2015)\nocite{wangTB2015} have proved a special case of this theorem when $C_j = 1$. The proof in the Appendix extends their result to the case $C_j > 1$.

\begin{corollary}\label{thm:PRBound}
If $x^*$ is $\epsilon$-optimal then $V^{PR}(C,0) \geq 0.5(1-\epsilon) V^{CDLP}$.
\end{corollary}
\proof{Proof.}
From Theorem \ref{thm:SingleLegCR} we know that
\[  \sum_{j=1}^L V_j^{PR}(C,0) \geq 0.5 \sum_{j=1}^L  \sum_{k=1}^K\Lambda_k \sum_{n \in N_j} \sum_{S \in \cal{O}_k : n \in S} x_k^*(S) P^k(n,S) r_n = \sum_{k=1}^K \Lambda_k \sum_{S \in \cal{S}} \sum_{n \in S} x_k^*(S) P^k(n,S) r_n.\]

Since $\sum_{j=1}^L V_j^{PR}(C,0)$ is the expected total reward of $PR$, and the right hand side is at least $(1-\epsilon)V^{CDLP}$, the expected reward of $PR$ is at least $0.5(1-\epsilon) V^{CDLP}$. \halmos
\endproof

Since $PR$ optimally manages the capacity allocation for each resource $l$, the expected future reward that $PR$ earns from each resource must dominate that of $FCFS$, which is stated in the following theorem.
\begin{theorem}\label{thm:PRDominance}
$V^{PR}(c,t) \geq V^{FCFS}(c,t).$
\end{theorem}
\proof{Proof.}
It is easy to check that $\cal{R}_t^{kl}(S,z) \geq \mathbf{1}(c_l >0) R_t^{kl}(S,z)$. The result follows after we combine this condition with (\ref{eq:DynamicsFCFS}) and (\ref{eq:DynamicsPR}).
\halmos
\endproof

\section{The Optimized Primal Routing Algorithm} \label{sec:OPR}
In this section, we propose a new algorithm called the Optimized Primal Routing ($OPR$) algorithm.  The advantage of $OPR$ is that it is much more practical than $PR$.  It does not perform random routing of customers.  Therefore, it exhibits much more stable behavior.  It also never offers an assortment that contains a product with $0$ inventory.  At the same time, it retains the performance guarantee of $PR$. 

For $OPR$, we only make one reasonable assumption on choice models, namely, that if we remove all products with non-positive rewards from an assortment, the expected reward of the assortment does not decrease.
\begin{assumption}\label{OPRAssumption}
For any assortment $S$ and any reward values $r_n$, $n=1,2,...,N$, let $H = \{ n \in S : r_n > 0\}$ be the set of products in $S$ with positive reward values. We must have for any customer type $k$,
\[ \sum_{n \in H} r_n P^k(n,H \cup \{0\}) \geq \sum_{n \in S} r_n P^k(n,S).\]
\end{assumption}
Note that this assumption holds for all random-utility models.

$OPR$ performs an additional optimization step compared to $PR$. Specifically, for each arriving customer, the algorithm locates an assortment that is at least as good as the assortment given by $PR$, according to the marginal values $\Delta_l V^{PR}$, $l=1,\ldots,L$, that are calculated in the same way as in $PR$:

\newcommand{\marginalallocation}{(Marginal Allocation Step) Let $c = (c_1,c_2,...,c_L)$ be the vector of inventory at time $t$. Upon an arrival of a type-$k$ customer at time $t$, offer an assortment $A_k^{OPR}(c,t)$ that aims at maximizing the marginal reward
\begin{equation}\label{eq:OPRStepTwo}
 A_k^{OPR}(c,t) \in \argmax_{S\in \cal{S}} \{ \sum_{l=1}^L R_t^{kl}(S, \Delta_l V^{PR}(c,t))\}.
 \end{equation}
When $S$ is a deterministic assortment the marginal reward becomes
 \[ \sum_{l=1}^L R_t^{kl}(S, \Delta_l V^{PR}(c,t)) = \lambda_k(t) \sum_{n \in S} P^k(n,S) \cdot (r_n - \Delta_l V^{PR}(c,t)).\]
Thus, \eqref{eq:OPRStepTwo} is equivalent to solving an assortment-optimization problem for a single customer with $r_n - \Delta_l V^{PR}(c,t)$ being the price for product $n$. In this step, $A_k^{OPR}(c,t)$ can be found by any approximation algorithms or heuristics. We only require that the expected marginal reward of  $A_k^{OPR}(c,t)$ be at least the marginal reward that $PR$ earns from this customer. Mathematically, we require that 
\begin{equation}\label{eq:OPRRequirement}
 \sum_{l=1}^L R_t^{kl}(A_k^{OPR}(c,t), \Delta_l V^{PR}(c,t)) \geq \sum_{l=1}^L \cal R_t^{kl}(A_k^{Static},  \Delta_l V^{PR}(c,t)).
\end{equation}
Note that it is trivial to satisfy this requirement because according to Assumption \ref{OPRAssumption}, we could just obtain $A_k^{OPR}(c,t)$ by taking the assortment in $\cal{O}_k$ that has the largest marginal reward, and then removing products with non-positive price $r_n - \Delta_l V^{PR}(c,t)$. The idea of \eqref{eq:OPRRequirement} is to enhance the empirical performance $OPR$ even further by conducting a broader search.
}

\begin{enumerate}
\item (Pre-Processing Step) Same as in FCFS.

\item \marginalallocation
\end{enumerate}

The requirement (\ref{eq:OPRRequirement}) also implies that $OPR$ satisfies the dynamic-substitution assumption. Therefore, the expected future reward $V^{OPR}(c,t)$ of $OPR$ satisfies the dynamic equation (\ref{eq:Dynamics}).

\begin{theorem} \label{thm:OPR}
The expected total reward of $OPR$ dominates that of $PR$.  That is, $V^{OPR}(C,0) \geq V^{PR}(C,0)$.
\end{theorem}
The proof defines a series of algorithms $\Pi^{(i)}$ that are intermediate to $PR$ and $OPR$.  Algorithm $\Pi^{(i)}$ is defined as follows.  For the first $i$ customers, apply \textit{OPR}. Afterward, for the $(i+1)$-th, $(i+2)$-th,..., customers, apply \textit{PR}. Thus, $\Pi^{(0)}$ resembles $OR$ and $\Pi^{(\infty)}$ resembles $OPR$. The proof in the Appendix shows that the expected revenue of each algorithm $\Pi^{(i)}$ improves upon that of $\Pi^{(i-1)}$.

Using Theorems \ref{thm:PRDominance} and \ref{thm:OPR} and Corollary \ref{thm:PRBound}, we can obtain our main result for $OPR$.
\begin{corollary}
If $x^*$ is $\epsilon$-optimal then $OPR$ is $0.5(1-\epsilon)$-competitive. Moreover, $OPR$ is asymptotically optimal as we scale up the total demand and capacity simultaneously.
\end{corollary}
\proof{Proof.}
The competitive ratio is a direct result of Theorem \ref{thm:OPR} and Corollary \ref{thm:PRBound}. It is well known that $FCFS$ is asymptotically optimal. Since $PR$ dominates $FCFS$ and $OPR$ dominates $PR$, $OPR$ must also be asymptotically optimal. \halmos
\endproof

\section{Appendix}

\noindent {\bf Proof of Theorem \ref{thm:upperbound}.}
Let $\delta_k$ be the actual number of arrivals of type-$k$ customers during the entire horizon. Let $I_k(S)$ be the number of times that a customer of type $k$ is shown assortment $S$ under \textit{OFF}. Let $J_{kn}$ be the number of times that a customer of type $k$ purchases product $n$ under \textit{OFF}.

We must have
\begin{equation}\label{eq:LPProof1}
 \sum_{S \in \cal{S}} I_k(S) = \delta_k,\,\,\, \forall k = 1,2,3...,K,
 \end{equation}
\begin{equation}\label{eq:LPProof2}
 \sum_{n=1}^N \sum_{k=1}^K J_{kn} \cdot \mathbf{1}(l_n = j) \leq C_j, \,\,\,\forall j = 1,2,...,M,
 \end{equation}
 \begin{equation}\label{eq:LPProof3}
 \mathbf{E}[J_{kn}] = \mathbf{E}[\sum_{S \in \cal{S} : n \in S} I_k(S) \cdot P^k(n,S)].
 \end{equation}
 
Note that equation (\ref{eq:LPProof3}) is a result of the dynamic-substitution assumption.
 
 Taking expectation on both sides of (\ref{eq:LPProof1}), we get
 \begin{equation}\label{eq:LPProof4}
  \sum_{S \in \cal{S}} \mathbf{E}[I_k(S)] = \Lambda_k,\,\,\, \forall k = 1,2,3...,K.
  \end{equation}
  
  Taking expectation on both sides of (\ref{eq:LPProof2}), we get
  \[ \sum_{n=1}^N \sum_{k=1}^K \mathbf{E}[J_{kn}] \cdot \mathbf{1}(l_n = j) \leq C_j \]
  \[ \Longrightarrow \sum_{n=1}^N \sum_{k=1}^K \mathbf{E}[\sum_{S \in \cal{S} : n \in S} I_k(S) \cdot P^k(n,S)] \cdot \mathbf{1}(l_n = j) \leq C_j \]
  \begin{equation}\label{eq:LPProof5}
  \Longrightarrow \sum_{S \in \cal{S}} \sum_{n \in S} \sum_{k=1}^K \mathbf{E}[ I_k(S) ] \cdot P^k(n,S)\mathbf{1}(l_n = j) \leq C_j .
  \end{equation}
  
  From (\ref{eq:LPProof4}) and (\ref{eq:LPProof5}) we know that $\mathbf{E}[ I_k(S) ]/ \Lambda_k$ is a feasible solution to LP (\ref{eq:ChoiceBasedLP}). Thus 
  \[ \mathbf{E}[R^{OFF}] =\sum_{k=1}^K \sum_{S\in \cal{S}} \Lambda_k \mathbf{E}[ I_k(S) ] \sum_{n\in S} P^k(n,S) r_n \]
  is at most the optimal objective value of (\ref{eq:ChoiceBasedLP}).
  \halmos

\noindent {\bf Proof of Theorem \ref{thm:nearOptimalCG}.}
When the column-generation algorithm terminates, we must have, $\forall k = 1,2,...,K$,
\begin{align*}
\begin{split}
 &\Lambda_k \sum_{n \in S_k^d, n>0} [r_n - \pi(l_n)] P^k(n,S_k^d)  - \sigma(k) \leq 0 ,\\
\Longrightarrow  &\Lambda_k \sum_{n \in S_k^*(\pi), n>0} [r_n - \pi(l_n)] P^k(n,S_k^*(\pi)) \cdot \frac{1}{1 + \epsilon}  - \sigma(k) \leq0\\
 \Longrightarrow &\Lambda_k \sum_{n \in S_k^*(\pi), n>0} [r_n - \pi(l_n)] P^k(n,S_k^*(\pi))  - \sigma(k) \leq \epsilon \sigma(k)\\
 \Longrightarrow &\Lambda_k \sum_{n \in S, n>0} [r_n - \pi(l_n)] P^k(n,S)  - \sigma(k) \leq \epsilon \sigma(k), \,\,\,\forall k = 1,2,...,K, \ \forall S \in \cal{S}.
 \end{split}
\end{align*}

Now let us look into the dual formulation of the CDLP given by  (\ref{eq:ChoiceBasedLP}):
\begin{align}
\begin{split} \label{eq:dualofChoiceBasedLP}
V^{CDLP}\equiv \min &\,\,\,\,\,\,\, \sum_j C_j\pi(j)+\sum_k \sigma(k) \\
\text{s.t. } &  \Lambda_k \sum_{n \in S, n>0} [r_n - \pi(l_n)] P^k(n,S)  - \sigma(k) \leq 0 , \,\,\,\forall k = 1,2,...,K, \ \forall S \in \cal{S},\\
 &\pi(j) \geq 0, \,\,\, \forall j=1,2,...,L\\
 & \sigma(k) \geq 0, \,\,\, \forall k =1,2,...,K.
 \end{split}
\end{align}
Let $\pi(j), j=1,...,L$ and $\sigma(k), k=1,...,K$, be optimal dual values for the reduced CDLP at termination of the column-generation algorithm.  Clearly, these variables will satisfy the dual constraint if we relax the first constraint in the linear program (\ref{eq:dualofChoiceBasedLP}) by 
\[  \Lambda_k \sum_{n \in S, n>0} [r_n - \pi(l_n)] P^k(n,S)  - \sigma(k) \leq \epsilon \sigma(k), \,\,\,\forall k = 1,2,...,K, \forall S \in \cal{S}.\]
That is, the optimal dual variables $\pi(j), j=1,...,L$ and $\sigma(k), k=1,...,K$ are feasible for the relaxed LP defined by
\begin{align}
\begin{split} \label{eq:reldualofChoiceBasedLP}
V^R\equiv  \min &\,\,\,\,\,\,\, \sum_j C_j\pi(j)+\sum_k \sigma(k) \\
\text{s.t. } & \Lambda_k \sum_{n \in S, n>0} [r_n - \pi(l_n)] P^k(n,S)  - \sigma(k) \leq \epsilon \sigma(k), \,\,\,\forall k = 1,2,...,K, \  \forall S \in \cal{S},\\
 &\pi(j) \geq 0, \,\,\, \forall j=1,2,...,L\\
 & \sigma(k) \geq 0, \,\,\, \forall k =1,2,...,K.
 \end{split}
\end{align}
Let us write down the primal of the relaxed dual program (\ref{eq:reldualofChoiceBasedLP}):
\begin{align}
\begin{split} \label{eq:dualrelaxChoiceBasedLP}
V^R\equiv \max & \,\,\,\,\,\,\, \sum_{k=1}^K \sum_{S\in \cal{S}} \Lambda_k x_k(S) \sum_{n\in S} P^k(n,S) r_n \\
\text{s.t. } &  \sum_{k =1}^K  \Lambda_k \sum_{n \in N_j} \sum_{S \in \cal{S}: n \in S} x_k(S) P^k(n,S) \leq C_j , \,\,\, \forall j = 1,2,...,L,\\
 &(1+\epsilon) \sum_{S \in \cal{S}} x_k(S) \leq 1, \,\,\, \forall k =1,2,...,K;\\
 & x_k(S) \geq 0, \,\,\, \forall k =1,2,...,K,\ \forall S \in \cal{S}.
 \end{split}
\end{align}
By strong duality, we know the optimal value returned by (\ref{eq:reldualofChoiceBasedLP}) and (\ref{eq:dualrelaxChoiceBasedLP}) should be the same and we denote it by $V^R$. Since (\ref{eq:reldualofChoiceBasedLP}) is less restrictive than (\ref{eq:dualofChoiceBasedLP}), then $V^R\leq V^{CDLP}$. Note that compared with the CDLP given in  (\ref{eq:ChoiceBasedLP}), (\ref{eq:dualrelaxChoiceBasedLP}) has nothing changed except that the left hand side of the second constraint is multiplied by $1+\epsilon$ and thus the entire problem becomes more restrictive. Now let us multiply the left hand side of the first constraint by $1+\epsilon$
\begin{align}
\begin{split} \label{eq:morerelaxedChoiceBasedLP}
 V^{R'}\equiv \max &\,\,\,\,\,\,\, \sum_{k=1}^K \sum_{S\in \cal{S}} \Lambda_k x_k(S) \sum_{n\in S} P^k(n,S) r_n \\
\text{s.t. } &  (1+\epsilon)\sum_{k =1}^K  \Lambda_k \sum_{n \in N_j}  \sum_{S \in \cal{S}: n \in S} x_k(S) P^k(n,S) \leq C_j , \,\,\, \forall j = 1,2,...,L,\\
 &(1+\epsilon) \sum_{S \in \cal{S}} x_k(S) \leq 1, \,\,\, \forall k =1,2,...,K;\\
 & x_k(S) \geq 0, \,\,\, \forall k =1,2,...,K,\ \forall S \in \cal{S}.
 \end{split}
\end{align}
Obviously, program (\ref{eq:morerelaxedChoiceBasedLP}) should have the same offering sets as those of (\ref{eq:ChoiceBasedLP}) at optimality, and the optimal value satisfies $V^{R'}=\frac{1}{1+\epsilon}V^{CDLP}$. Moreover, since (\ref{eq:morerelaxedChoiceBasedLP}) is more restrictive than (\ref{eq:dualrelaxChoiceBasedLP}), we know $V^R \geq V^{R'}$. Given that the values of current $\pi(j), j=1,...,L$ and $\sigma(k), k=1,...,K$ are feasible to (\ref{eq:reldualofChoiceBasedLP}), we have $\sum_j C_j\pi(j)+\sum_k \sigma(k) \geq V^R$. Combining all the information together, we know $\sum_j C_j\pi(j)+\sum_k \sigma(k)>\frac{1}{1+\epsilon}V^{CDLP}>(1-\epsilon)V^{CDLP}$. In other words, the value obtained from the column- generation algorithm $\sum_j C_j\pi(j)+\sum_k \sigma(k)$ is at least $1-\epsilon$ times $V^{CDLP}$.
\halmos

\noindent {\bf Proof of Theorem \ref{thm:FCFS}.}

Let \[ s_{kn}^* \equiv \sum_{S \in \cal{O}_k : n \in S} x_k^*(S) P^k(n,S)\]be the probability that a customer of type $k$ chooses product $n$ from $A_k^{Static}$.

Under $FCFS$, the total number of customers who will choose resource $j$, including those who are rejected due to a lack of inventory, is a Poisson random variable $D_j$ with mean
\[ \mathbf{E}[D_j] = \sum_{k=1}^K\Lambda_k \sum_{n\in N_j} s_{kn}^*,\]
which is at most $C_j$ according to the capacity constraint of (\ref{eq:ChoiceBasedLP}).

Conditioned on the event that a customer of type $k$ successfully purchases a product associated with resource $j$, the probability that the purchased product is $i$ is
\[ \frac{ \mathbf{1}(l_i = j) s_{ki}^*}{ \sum_{n\in N_j}  s_{kn}^*},\]
which is independent of the choice of other customers.

Then, conditioned on the event that a customer of type $k$ successfully purchases a product associated with resource $j$, the expected reward that the customer brings is
\begin{equation}\label{eq:proofFCFS1} \frac{\sum_{n \in N_j} s_{kn}^* r_n}{ \sum_{n \in N_j} s_{kn}^*}.\end{equation}

According to the properties of Poisson processes, conditioned on the event $D_j = d$, each of the $d$ customers can be seen as randomly and independently picked from the entire horizon. The probability that each of the $d$ customers is of type $k$ is
\begin{equation}\label{eq:proofFCFS2}  \frac{\Lambda_k \sum_{n \in N_j} s_{kn}^*}{ \sum_{i=1}^K \Lambda_i \sum_{n \in N_j} s_{in}^*} =  \frac{\Lambda_k \sum_{n \in N_j} s_{kn}^*}{\mathbf{E}[D_j]}.\end{equation}

Under FCFS, for any integer $d \leq C_j$, if $D_j = d$, all of these $d$ customer demands will be satisfied. Combining (\ref{eq:proofFCFS1}) and (\ref{eq:proofFCFS2}), we know that for any $d \leq C_j$,
\begin{align*}
& \mathbf{E}[\text{Total reward obtained from resource $j$} | D_j = d ] \\
= & d \cdot \mathbf{E}[\text{Total reward obtained from resource $j$} | D_j = 1]\\
= & d \cdot \sum_{k=1}^K  \frac{\sum_{n \in N_j} s_{kn}^* r_n}{ \sum_{n \in N_j} s_{kn}^*} \cdot \frac{\Lambda_k \sum_{n \in N_j} s_{kn}^*}{\mathbf{E}[D_j]}\\
= & d \cdot  \frac{\sum_{k=1}^K  \Lambda_k \sum_{n \in N_j} s_{kn}^* r_n}{\mathbf{E}[D_j]}.
\end{align*}

Then,
\begin{align*}
&\mathbf{E}[\text{Total reward obtained from resource $j$}]\\
= &\sum_{d=0}^{\infty} P(D_j = d)  \times \mathbf{E}[\text{Total reward obtained from resource }j |D_j = d]\\
\geq &\sum_{d=0}^{C_j} P(D_j = d)  \times \mathbf{E}[\text{Total reward obtained from resource }j |D_j = d]\\
 = & \sum_{d=0}^{C_j} \frac{\mathbf{E}[D_j]^d}{d!}e^{-\mathbf{E}[D_j]} \times d \cdot  \frac{ \sum_{k=1}^K\Lambda_k \sum_{n \in N_j} s_{kn}^* r_n}{ \mathbf{E}[D_j]}\\
 \geq & \frac{1}{e} \cdot \sum_{k=1}^K\Lambda_k \sum_{n \in N_j} s_{kn}^* r_n.
\end{align*}
The last step follows from the fact that for any non-negative value $x$, it always holds that
\[ \sum_{i=0}^{\lceil x \rceil} \frac{x^i}{i!} e^{-x}\cdot \frac{i}{x} \geq \frac{1}{e}.\]

Thus, the expected total reward that FCFS earns from all resources is at least
\[ \sum_{j=1}^L \frac{1}{e} \cdot \sum_{k=1}^K\Lambda_k \sum_{n \in N_j} s_{kn}^* r_n = \frac{1}{e} \cdot \sum_{k=1}^K\Lambda_k \sum_{n=1}^N  s_{kn}^* r_n,\]
which is $1/e$ times the objective value of (\ref{eq:ChoiceBasedLP}) corresponding to solution $x_k^*(S)$. Since $x_k^*(S)$ is $\epsilon$-optimal, and $V^{CDLP}$ is an upper bound on $V^{OFF}$, the total expected reward that FCFS earns is at least $\frac{1}{e}(1-\epsilon)\cdot V^{OFF}$.
\endproof

\noindent {\bf Proof of Theorem \ref{thm:SingleLegCR}.}
Wang, Truong and Bank (2015)\nocite{wangTB2015} have proved a special case of this theorem when $C_j = 1$. In this proof we extend their results to the case $C_j > 1$.

Let \[ s_{kn}^* \equiv \sum_{S \in \cal{O}_k : n \in S} x_k^*(S) P^k(n,S)\]be the probability that a customer of type $k$ chooses product $n$ from $A_k^{Static}$.

The total expected demand $\mathbf{E}[D_j]$ that is allocated to resource $j$ is defined by the optimal solution $x^*$ of (\ref{eq:ChoiceBasedLP})
\[ \mathbf{E}[D_j] \equiv \sum_{k=1}^K \Lambda_k \sum_{n \in N_j} \sum_{S \in \cal{O}_k : n \in S}  x_k^*(S) P^k(n,S) = \sum_{k=1}^K\Lambda_k \sum_{n \in N_j} s_{kn}^*,\]
which is at most $C_j$ according to the constraints of (\ref{eq:ChoiceBasedLP}).

Consider the following sub-optimal policy applied to the single-resource revenue-management problem for resource $j$. Divide the horizon $[0,1]$ into $C_j$ intervals $[t(0),t(1)]$, $[t(1), t(2)]$,...,$[t(C_j-1),t(C_j)]$, such that $t(0) = 0$, $t(C_j) = 1$ and 
\begin{equation}\label{eq:PRProofa}
 \int_{t(i-1)}^{t(i)}  \sum_{k=1}^K \lambda_k(u) \sum_{n\in N_j} s_{kn}^*  du = \frac{\mathbf{E}[D_j]}{C_j} \leq 1, \,\, \forall i = 1,2,...,C_j. 
 \end{equation}
In other words, the average number of arrivals at resource $j$ in each interval has the same value $\mathbf{E}[D_j] / C_j$.   Then, we view the $C_j$ units of resource $j$ as $C_j$ different resources with unit inventory. The $i$-th resource only accepts customers, if any, arriving during the $i$-th interval $[t(i-1),t(i)]$, for $i=1,2,...,C_j$. We optimally solve the admission-control problem for each resource with unit inventory. Let $g_i(t)$ be the expected future reward obtained from the $i$-th resource. It satisfies
\begin{equation}\label{eq:PRProofb}
\frac{dg_i(t)}{dt} = -  \sum_{k=1}^K \cal R_t^{kj}(A_k^{Static},g_i(t)), \,\,\, \forall t \in (t(i-1),t(i))
 \end{equation}
with boundary condition
\begin{equation}\label{eq:PRProofc}
 g_i(t(i)) = 0.
 \end{equation}

Wang, Truong and Bank (2015)\nocite{wangTB2015} have shown that, once conditions (\ref{eq:PRProofa}), (\ref{eq:PRProofb}) and (\ref{eq:PRProofc}) hold at the same time, the total expected reward obtained from resource $i$ is at least
\[ g_i(t(i-1)) \geq 0.5 \int_{t(i-1)}^{t(i)}  \sum_{k=1}^K \lambda_k(u) \sum_{n \in N_j} s_{kn}^* r_n  du.\]

Thus, under this sub-optimal policy, the expected total reward obtained from all $C_j$ resources is at least
\[ \sum_{i=1}^{C_j} g_i(t(i-1)) \geq 0.5 \sum_{i=1}^{C_j} \int_{t(i-1)}^{t(i)}  \sum_{k=1}^K \lambda_k(u) \sum_{n \in N_j} s_{kn}^* r_n  du = 0.5 \sum_{k=1}^K \Lambda_k \sum_{n \in N_j} s_{kn}^* r_n .\]
Since $V_j^{PR}(C,0)$ is the optimal expected reward when the decisions are made for the entire resource $j$, $V_j^{PR}(C,0)$ must be at least the expected total reward of the sub-optimal policy. In other words.
\[  V_j^{PR}(C,0) \geq 0.5 \sum_{k=1}^K \Lambda_k \sum_{n \in N_j} s_{kn}^* r_n.\]

\halmos

\noindent {\bf Proof of Theorem \ref{thm:OPR}.} 

We want to show that
\[ V^{OPR}(c,t) \geq V^{PR}(c,t)\]
for every given state $(c,t)$.

 Define an algorithm $\Pi^{(i)}$ as follows. For the first $i$ customers, apply \textit{OPR}. Afterward, for the $(i+1)$-th, $(i+2)$-th,..., customers, apply \textit{PR}. Let $h^{(i)}(c,t)$ be the expected future reward when policy $\Pi^{(i)}$ is applied starting at time $t$ with remaining inventory $c(t)$, and \textit{assuming that no customers have arrived prior to time $t$}. We must have
\[ h^{(0)}(c,t) = V^{PR}(c,t),\]
\[ \lim_{i\to \infty} h^{(i)}(c,t) = V^{OPR}(c,t) .\]

The dynamic programming equation for algorithm $\Pi^{(1)}$ is
\begin{align}\label{eq:OptimizedAlgProof2}
\begin{split}
\frac{\partial h^{(1)}(c,t)}{\partial t} &= - \sum_{k=1}^K \left[ \sum_{l=1}^L  R_t^{kl}(A_k^{OPR}(c,t), h^{(1)}(c,t) - h^{(0)}(c-e_l,t)) - \lambda_k(t) \mathbf{E}[P^k(0,A_k^{OPR}(c,t))] \Delta^{(1)}(c,t) \right] \\
& = - \sum_{k=1}^K \left[ \sum_{l=1}^L  R_t^{kl}(A_k^{OPR}(c,t), \Delta_l V^{PR}(c,t)) - \lambda_k(t) \Delta^{(1)}(c,t) \right],
\end{split}
\end{align}
 where
\[ \Delta^{(1)}(c,t) \equiv h^{(1)}(c,t) - V^{PR}(c,t)\]
is the difference in the expected future reward between $\Pi^{(1)}$ and \textit{PR} (\textit{PR}  $\equiv \Pi^{(0)}$).

Combining (\ref{eq:OPRRequirement}), (\ref{eq:OptimizedAlgProof2}) and the dynamic equation (\ref{eq:DynamicsPR}) for $PR$, we can obtain
\begin{equation}\frac{\partial h^{(1)}(c,t)}{\partial t} \leq  \frac{\partial V^{PR}(c,t)}{\partial t} + \sum_{k=1}^K \lambda_k(t)  \Delta^{(1)}(c,t).  \end{equation}
This equation implies that, if at some time $t_0$ we have $\Delta^{(1)}(c,t_0)< 0$ or equivalently
\begin{equation} \label{eq:OptimizedAlgProofb}
h^{(1)}(c,t_0) - V^{PR}(c,t_0) < 0,
 \end{equation}
then we must have
\begin{equation}  \frac{\partial h^{(1)}(c,t)}{\partial t} <  \frac{\partial V^{PR}(c,t)}{\partial t}, \,\,\, \forall t \in (t_0,1] \end{equation}
and
\begin{equation} \label{eq:OptimizedAlgProofa}
 h^{(1)}(c,t) < V^{PR}(c,t), \,\,\, \forall t \in (t_0,1].
 \end{equation}
However, since we know that $h^{(1)}(c,1) = V^{PR}(c,1) = 0$, (\ref{eq:OptimizedAlgProofa}) cannot be true, and thus (\ref{eq:OptimizedAlgProofb}) cannot be true. Therefore, we have proved 

\begin{equation}\label{eq:OptimizedAlgProof4}
h^{(1)}(c,t) \geq V^{PR}(c,t), \,\,\, \forall t \in [0,1].
 \end{equation}
 
Next, we show that 
\begin{equation}\label{eq:OptimizedAlgProof5}
h^{(i)}(c,t) \geq h^{(i-1)}(c,t), \,\,\, \forall t\in [0,1]
\end{equation}
 by induction on $i$. 
 
 Equation (\ref{eq:OptimizedAlgProof4}) already proves the base case $i=1$. Suppose for some $\bar i > 1$, (\ref{eq:OptimizedAlgProof5}) holds for all $i < \bar i$. Now we show that it also holds for $i = \bar i$. By definition, for any $\bar i > 1$, algorithms $\Pi^{(\bar i)}$ and $\Pi^{(\bar i-1)}$ must offer the same assortment to the first customer, for they both apply \textit{OPR} to the first customer. Thus, $\Pi^{(\bar i)}$ and $\Pi^{(\bar i-1)}$ earn the same reward from the first customer, and then transit into the same state. After that first customer,  $\Pi^{(\bar i)}$ continues to apply $\Pi^{(\bar i-1)}$ pretending that no customer has ever arrived, while $\Pi^{(\bar i-1)}$ continues to apply $\Pi^{(\bar i-2)}$. By induction, the expected future reward of  $\Pi^{(\bar i-1)}$ is at least that of $\Pi^{(\bar i-2)}$. Therefore, the expected future reward of $\Pi^{(\bar i)}$ is at least that of $\Pi^{(\bar i-1)}$.

Thus, we have proved (\ref{eq:OptimizedAlgProof5}). It immediately follows that 
\[ h^{(\infty)}(c,t) \geq h^{(0)}(c,t) \Longrightarrow V^{OPR}(c,t) \geq V^{PR}(c,t).\]
\halmos

\bibliographystyle{ormsv080}
\bibliography{myrefs}{}

\begin{thebibliography}{57}
\expandafter\ifx\csname natexlab\endcsname\relax\def\natexlab#1{#1}\fi
\expandafter\ifx\csname url\endcsname\relax
  \def\url#1{{\tt #1}}\fi
\expandafter\ifx\csname urlprefix\endcsname\relax\def\urlprefix{URL }\fi
\expandafter\ifx\csname urlstyle\endcsname\relax
  \expandafter\ifx\csname doi\endcsname\relax
  \def\doi#1{doi:\discretionary{}{}{}#1}\fi \else
  \expandafter\ifx\csname doi\endcsname\relax
  \def\doi{doi:\discretionary{}{}{}\begingroup \urlstyle{rm}\Url}\fi \fi

\bibitem[{Agrawal et~al.(2009)Agrawal, Wang, and Ye}]{agrawal2009dynamic}
Agrawal, Shipra, Zizhuo Wang, Yinyu Ye. 2009.
\newblock A dynamic near-optimal algorithm for online linear programming.
\newblock {\it arXiv preprint arXiv:0911.2974\/} .

\bibitem[{Babaioff et~al.(2008)Babaioff, Immorlica, Kempe, and
  Kleinberg}]{babaioff2008online}
Babaioff, Moshe, Nicole Immorlica, David Kempe, Robert Kleinberg. 2008.
\newblock Online auctions and generalized secretary problems.
\newblock {\it ACM SIGecom Exchanges\/} {\bf 7}(2) 7.

\bibitem[{Bahmani and Kapralov(2010)}]{bahmani2010improved}
Bahmani, Bahman, Michael Kapralov. 2010.
\newblock Improved bounds for online stochastic matching.
\newblock {\it Algorithms--ESA 2010\/}. Springer, 170--181.

\bibitem[{Ball and Queyranne(2009)}]{ball2009toward}
Ball, Michael~O, Maurice Queyranne. 2009.
\newblock Toward robust revenue management: Competitive analysis of online
  booking.
\newblock {\it Operations Research\/} {\bf 57}(4) 950--963.

\bibitem[{Belobaba(1989)}]{belobaba1989}
Belobaba, Peter~P. 1989.
\newblock Or practice—application of a probabilistic decision model to
  airline seat inventory control.
\newblock {\it Operations Research\/} {\bf 37}(2) 183--197.
\newblock \doi{10.1287/opre.37.2.183}.
\newblock \urlprefix\url{http://dx.doi.org/10.1287/opre.37.2.183}.

\bibitem[{Bernstein et~al.(2010)Bernstein, Kok, and
  Xie}]{Bernstein10dynamicassortment}
Bernstein, Fernando, A.~Gurhan Kok, Lei Xie. 2010.
\newblock Dynamic assortment customization with limited inventories.

\bibitem[{Blanchet et~al.(2013)Blanchet, Gallego, and Goyal}]{Blanchet:2013}
Blanchet, Jose, Guillermo Gallego, Vineet Goyal. 2013.
\newblock A markov chain approximation to choice modeling.
\newblock {\it Proceedings of the Fourteenth ACM Conference on Electronic
  Commerce\/}. EC '13, ACM, New York, NY, USA, 103--104.
\newblock \doi{10.1145/2482540.2482560}.
\newblock \urlprefix\url{http://doi.acm.org/10.1145/2482540.2482560}.

\bibitem[{Bront et~al.(2009)Bront, Mendez-Diaz, and Vulcano}]{BrontMV2009}
Bront, J.J.M., I.~Mendez-Diaz, G.~Vulcano. 2009.
\newblock {A column generation algorithm for choice-based network revenue
  management}.
\newblock {\it Operations Research\/} {\bf 57}(3) 769--784.

\bibitem[{Cardoen et~al.(2010)Cardoen, Demeulemeester, and
  Beli{\"e}n}]{cardoen2010operating}
Cardoen, Brecht, Erik Demeulemeester, Jeroen Beli{\"e}n. 2010.
\newblock Operating room planning and scheduling: A literature review.
\newblock {\it European Journal of Operational Research\/} {\bf 201}(3)
  921--932.

\bibitem[{Chan and Farias(2009)}]{chan2009stochastic}
Chan, Carri~W, Vivek~F Farias. 2009.
\newblock Stochastic depletion problems: Effective myopic policies for a class
  of dynamic optimization problems.
\newblock {\it Mathematics of Operations Research\/} {\bf 34}(2) 333--350.

\bibitem[{Dai et~al.(2014)Dai, Ding, Kleywegt, Wang, and Zhang}]{JimDai2014}
Dai, Jim, Weijun Ding, Anton~J. Kleywegt, Xinchang Wang, Yi~Zhang. 2014.
\newblock Choice based revenue management for parallel flights.
\newblock Tech. rep., Working Paper.

\bibitem[{Devanur(2009)}]{Devanur09theadwords}
Devanur, Nikhil~R. 2009.
\newblock The adwords problem: Online keyword matching with budgeted bidders
  under random permutations.
\newblock {\it In Proc. 10th Annual ACM Conference on Electronic Commerge
  (EC\/}.

\bibitem[{Devanur et~al.(2011)Devanur, Jain, Sivan, and
  Wilkens}]{devanur2011near}
Devanur, Nikhil~R, Kamal Jain, Balasubramanian Sivan, Christopher~A Wilkens.
  2011.
\newblock Near optimal online algorithms and fast approximation algorithms for
  resource allocation problems.
\newblock {\it Proceedings of the 12th ACM conference on Electronic
  commerce\/}. ACM, 29--38.

\bibitem[{Duhigg(2012)}]{duhigg2012companies}
Duhigg, Charles. 2012.
\newblock How companies learn your secrets.
\newblock {\it The New York Times\/} {\bf 16} 2012.

\bibitem[{Feldman et~al.(2014)Feldman, Liu, Topaloglu, and
  Ziya}]{feldman2014appointment}
Feldman, Jacob, Nan Liu, Huseyin Topaloglu, Serhan Ziya. 2014.
\newblock Appointment scheduling under patient preference and no-show behavior.
\newblock {\it Operations Research\/} {\bf 62}(4) 794--811.

\bibitem[{Feldman et~al.(2010)Feldman, Henzinger, Korula, Mirrokni, and
  Stein}]{feldman2010online}
Feldman, Jon, Monika Henzinger, Nitish Korula, Vahab~S Mirrokni, Cliff Stein.
  2010.
\newblock Online stochastic packing applied to display ad allocation.
\newblock {\it Algorithms--ESA 2010\/}. Springer, 182--194.

\bibitem[{Feldman et~al.(2009)Feldman, Mehta, Mirrokni, and
  Muthukrishnan}]{feldman2009online}
Feldman, Jon, Aranyak Mehta, Vahab Mirrokni, S~Muthukrishnan. 2009.
\newblock Online stochastic matching: Beating 1-1/e.
\newblock {\it Foundations of Computer Science, 2009. FOCS'09. 50th Annual IEEE
  Symposium on\/}. IEEE, 117--126.

\bibitem[{Gallego et~al.(2004)Gallego, Iyengar, Phillips, and
  Dubey}]{gallego2004managing}
Gallego, Guillermo, Garud Iyengar, R~Phillips, Abha Dubey. 2004.
\newblock Managing flexible products on a network.
\newblock Unpublished.

\bibitem[{Gallego et~al.(2014)Gallego, Ratliff, and
  Shebalov}]{gallego2014general}
Gallego, Guillermo, Richard Ratliff, Sergey Shebalov. 2014.
\newblock A general attraction model and sales-based linear program for network
  revenue management under customer choice.
\newblock {\it Operations Research\/} .

\bibitem[{Goel and Mehta(2008)}]{goel2008online}
Goel, Gagan, Aranyak Mehta. 2008.
\newblock Online budgeted matching in random input models with applications to
  adwords.
\newblock {\it Proceedings of the nineteenth annual ACM-SIAM symposium on
  Discrete algorithms\/}. Society for Industrial and Applied Mathematics,
  982--991.

\bibitem[{Golrezaei et~al.(2014)Golrezaei, Nazerzadeh, and
  Rusmevichientong}]{golrezaei2014real}
Golrezaei, Negin, Hamid Nazerzadeh, Paat Rusmevichientong. 2014.
\newblock Real-time optimization of personalized assortments.
\newblock {\it Management Science\/} {\bf 60}(6) 1532--1551.

\bibitem[{Guerriero and Guido(2011)}]{guerriero2011operational}
Guerriero, Francesca, Rosita Guido. 2011.
\newblock Operational research in the management of the operating theatre: a
  survey.
\newblock {\it Health care management science\/} {\bf 14}(1) 89--114.

\bibitem[{Gupta(2007)}]{gupta2007surgical}
Gupta, D. 2007.
\newblock {Surgical suites' operations management}.
\newblock {\it Production and Operations Management\/} {\bf 16}(6) 689--700.

\bibitem[{Gupta and Wang(2008)}]{gupta2008revenue}
Gupta, Diwakar, Lei Wang. 2008.
\newblock Revenue management for a primary-care clinic in the presence of
  patient choice.
\newblock {\it Operations Research\/} {\bf 56}(3) 576--592.

\bibitem[{Haeupler et~al.(2011)Haeupler, Mirrokni, and
  Zadimoghaddam}]{haeupler2011online}
Haeupler, Bernhard, Vahab~S Mirrokni, Morteza Zadimoghaddam. 2011.
\newblock Online stochastic weighted matching: Improved approximation
  algorithms.
\newblock {\it Internet and Network Economics\/}. Springer, 170--181.

\bibitem[{Honhon et~al.(2010)Honhon, Gaur, and Seshadri}]{honhonor2010}
Honhon, Dorothee, Vishal Gaur, Sridhar Seshadri. 2010.
\newblock Assortment planning and inventory decisions under stockout-based
  substitution.
\newblock {\it Operations Research\/} {\bf 58}(5) 1364--1379.

\bibitem[{Jaillet and Lu(2011)}]{jaillet2011online}
Jaillet, Patrick, Xin Lu. 2011.
\newblock Online resource allocation problems.
\newblock Tech. rep., Working paper.

\bibitem[{Jaillet and Lu(2013)}]{jaillet2013online}
Jaillet, Patrick, Xin Lu. 2013.
\newblock Online stochastic matching: New algorithms with better bounds.
\newblock {\it Mathematics of Operations Research\/} {\bf 39}(3) 624--646.

\bibitem[{Jasin and Kumar(2012)}]{JasinKumar2012}
Jasin, Stefanus, Sunil Kumar. 2012.
\newblock A re-solving heuristic with bounded revenue loss for network revenue
  management with customer choice.
\newblock {\it Mathematics of Operations Research\/} {\bf 37}(2) 313--345.
\newblock \doi{10.1287/moor.1120.0537}.
\newblock \urlprefix\url{http://dx.doi.org/10.1287/moor.1120.0537}.

\bibitem[{Karande et~al.(2011)Karande, Mehta, and Tripathi}]{karande2011online}
Karande, Chinmay, Aranyak Mehta, Pushkar Tripathi. 2011.
\newblock Online bipartite matching with unknown distributions.
\newblock {\it Proceedings of the forty-third annual ACM symposium on Theory of
  computing\/}. ACM, 587--596.

\bibitem[{Karp et~al.(1990)Karp, Vazirani, and Vazirani}]{karp1990optimal}
Karp, R.~M., U.~V. Vazirani, V.~V. Vazirani. 1990.
\newblock An optimal algorithm for on-line bipartite matching.
\newblock {\it Proceedings of the Twenty-second Annual ACM Symposium on Theory
  of Computing\/}. STOC '90, ACM, New York, NY, USA, 352--358.
\newblock \doi{10.1145/100216.100262}.
\newblock \urlprefix\url{http://doi.acm.org/10.1145/100216.100262}.

\bibitem[{Kleinberg(2005)}]{kleinberg2005multiple}
Kleinberg, Robert. 2005.
\newblock A multiple-choice secretary algorithm with applications to online
  auctions.
\newblock {\it Proceedings of the sixteenth annual ACM-SIAM symposium on
  Discrete algorithms\/}. Society for Industrial and Applied Mathematics,
  630--631.

\bibitem[{Kok and Fisher(2007)}]{kokfisher2007}
Kok, A.~Gurhan, Marshall~L. Fisher. 2007.
\newblock Demand estimation and assortment optimization under substitution:
  Methodology and application.
\newblock {\it Operations Research\/} {\bf 55}(6) 1001--1021.

\bibitem[{Kunnumkal and Topaloglu(2010)}]{KunnumkalT2010}
Kunnumkal, Sumit, Huseyin Topaloglu. 2010.
\newblock A new dynamic programming decomposition method for the network
  revenue management problem with customer choice behavior.
\newblock {\it Production and Operations Management\/} {\bf 19}(5) 575--590.
\newblock \doi{10.1111/j.1937-5956.2009.01118.x}.
\newblock \urlprefix\url{http://dx.doi.org/10.1111/j.1937-5956.2009.01118.x}.

\bibitem[{Lee and Hersh(1993)}]{lee1993model}
Lee, Tak~C, Marvin Hersh. 1993.
\newblock A model for dynamic airline seat inventory control with multiple seat
  bookings.
\newblock {\it Transportation Science\/} {\bf 27}(3) 252--265.

\bibitem[{Li et~al.(2015)Li, Rusmevichientong, and
  Topaloglu}]{doi:10.1287/opre.2015.1355}
Li, Guang, Paat Rusmevichientong, Huseyin Topaloglu. 2015.
\newblock The d-level nested logit model: Assortment and price optimization
  problems.
\newblock {\it Operations Research\/} {\bf 63}(2) 325--342.
\newblock \doi{10.1287/opre.2015.1355}.
\newblock \urlprefix\url{http://dx.doi.org/10.1287/opre.2015.1355}.

\bibitem[{Littlewood(2005)}]{littlewood2005special}
Littlewood, Ken. 2005.
\newblock Special issue papers: Forecasting and control of passenger bookings.
\newblock {\it Journal of Revenue and Pricing Management\/} {\bf 4}(2)
  111--123.

\bibitem[{Liu and van Ryzin(2008)}]{LiuV2008}
Liu, Q., G.~van Ryzin. 2008.
\newblock {On the choice-based linear programming model for network revenue
  management}.
\newblock {\it Manufacturing \& Service Operations Management\/} {\bf 10}(2)
  288.

\bibitem[{Mahajan and van Ryzin(2001)}]{MahajanV2001}
Mahajan, S., G.~van Ryzin. 2001.
\newblock {Stocking retail assortments under dynamic consumer substitution}.
\newblock {\it Operations Research\/} {\bf 49}(3) 334--351.

\bibitem[{Mahdian and Yan(2011)}]{mahdian2011online}
Mahdian, Mohammad, Qiqi Yan. 2011.
\newblock Online bipartite matching with random arrivals: an approach based on
  strongly factor-revealing lps.
\newblock {\it Proceedings of the forty-third annual ACM symposium on Theory of
  computing\/}. ACM, 597--606.

\bibitem[{Manshadi et~al.(2012)Manshadi, Gharan, and
  Saberi}]{manshadi2012online}
Manshadi, Vahideh~H, Shayan~Oveis Gharan, Amin Saberi. 2012.
\newblock Online stochastic matching: Online actions based on offline
  statistics.
\newblock {\it Mathematics of Operations Research\/} {\bf 37}(4) 559--573.

\bibitem[{May et~al.(2010)May, Spangler, Strum, and Vargas}]{may2010surgical}
May, J.H., W.E. Spangler, D.P. Strum, L.G. Vargas. 2010.
\newblock {The surgical scheduling problem: Current research and future
  opportunities}.
\newblock {\it Production and Operations Management\/} .

\bibitem[{Mayer(1976)}]{mayer1976}
Mayer, M. 1976.
\newblock {Seat allocation, or simple model of seat allocation via
  sophisticated ones}.
\newblock {\it 16th AGIFORS Symposium Proceedings\/}. 103--135.

\bibitem[{McFadden and Train(2000)}]{mcfadden2000}
McFadden, Daniel, Kenneth Train. 2000.
\newblock {Mixed MNL models for discrete response}.
\newblock {\it Journal of Applied Econometrics\/} {\bf 15}(5) 447--470.
\newblock
  \urlprefix\url{https://ideas.repec.org/a/jae/japmet/v15y2000i5p447-470.html}.

\bibitem[{Mehta(2012)}]{mehta2012online}
Mehta, Aranyak. 2012.
\newblock Online matching and ad allocation.
\newblock {\it Theoretical Computer Science\/} {\bf 8}(4) 265--368.

\bibitem[{Meissner and Strauss(2012)}]{Meissner2012459}
Meissner, Joern, Arne Strauss. 2012.
\newblock Network revenue management with inventory-sensitive bid prices and
  customer choice.
\newblock {\it European Journal of Operational Research\/} {\bf 216}(2) 459 --
  468.
\newblock \doi{http://dx.doi.org/10.1016/j.ejor.2011.06.033}.
\newblock
  \urlprefix\url{http://www.sciencedirect.com/science/article/pii/S0377221711005625}.

\bibitem[{Mittal and Schulz(2013)}]{doi:10.1287/opre.1120.1093}
Mittal, Shashi, Andreas~S. Schulz. 2013.
\newblock A general framework for designing approximation schemes for
  combinatorial optimization problems with many objectives combined into one.
\newblock {\it Operations Research\/} {\bf 61}(2) 386--397.
\newblock \doi{10.1287/opre.1120.1093}.
\newblock \urlprefix\url{http://dx.doi.org/10.1287/opre.1120.1093}.

\bibitem[{Molinaro and Ravi(2013)}]{molinaro2013geometry}
Molinaro, Marco, R~Ravi. 2013.
\newblock The geometry of online packing linear programs.
\newblock {\it Mathematics of Operations Research\/} {\bf 39}(1) 46--59.

\bibitem[{P.~Cachon et~al.(2005)P.~Cachon, Terwiesch, and Xu}]{P.Cachon2005}
P.~Cachon, G{\'e}rard, Christian Terwiesch, Yi~Xu. 2005.
\newblock Retail assortment planning in the presence of consumer search.
\newblock {\it Manufacturing \& Service Operations Management\/} {\bf 7}(4)
  330--346.
\newblock \doi{10.1287/msom.1050.0088}.
\newblock \urlprefix\url{http://dx.doi.org/10.1287/msom.1050.0088}.

\bibitem[{Qin et~al.(2015)Qin, Zhang, Hua, and Shi}]{CongApproximationRM}
Qin, Chao, Huanan Zhang, Cheng Hua, Cong Shi. 2015.
\newblock A simple admission control policy for revenue management problems
  with non-stationary customer arrivals.
\newblock {\it working paper\/} .

\bibitem[{Talluri and Van~Ryzin(2004)}]{TalluriV2004}
Talluri, K., G.~Van~Ryzin. 2004.
\newblock {Revenue management under a general discrete choice model of consumer
  behavior}.
\newblock {\it Management Science\/} {\bf 50}(1) 15--33.

\bibitem[{Talluri and Van~Ryzin(2005)}]{talluri2005theory}
Talluri, K.T., G.~Van~Ryzin. 2005.
\newblock {\it The theory and practice of revenue management\/}, vol.~68.
\newblock Springer Verlag, New York.

\bibitem[{Topaloglu(2013)}]{Topaloglu2013stock}
Topaloglu, Huseyin. 2013.
\newblock Joint stocking and product offer decisions under the multinomial
  logit model.
\newblock {\it Production and Operations Management\/} {\bf 22}(5) 1182--1199.
\newblock \doi{10.1111/j.1937-5956.2012.01423.x}.
\newblock \urlprefix\url{http://dx.doi.org/10.1111/j.1937-5956.2012.01423.x}.

\bibitem[{Van~Ryzin and Mahajan(1999)}]{VanRyzinM1999}
Van~Ryzin, G., S.~Mahajan. 1999.
\newblock {On the relationship between inventory costs and variety benefits in
  retail assortments}.
\newblock {\it Management Science\/} {\bf 45}(11) 1496--1509.

\bibitem[{Wang et~al.(2015)Wang, Truong, and Bank}]{wangTB2015}
Wang, X., V.~Truong, D.~Bank. 2015.
\newblock Online advance admission scheduling for services, with customer
  preferences.
\newblock Working paper.

\bibitem[{Zhang and Adelman(2009)}]{zhang2009af}
Zhang, Dan, Daniel Adelman. 2009.
\newblock An approximate dynamic programming approach to network revenue
  management with customer choice.
\newblock {\it Transportation Science\/} {\bf 43}(3) 381--394.
\newblock \doi{10.1287/trsc.1090.0262}.
\newblock \urlprefix\url{http://dx.doi.org/10.1287/trsc.1090.0262}.

\bibitem[{Zhang and Cooper(2005)}]{zhang2005revenue}
Zhang, Dan, William~L Cooper. 2005.
\newblock Revenue management for parallel flights with customer-choice
  behavior.
\newblock {\it Operations Research\/} {\bf 53}(3) 415--431.

\end{thebibliography}

\end{document}